\documentclass[11pt]{article}
\usepackage{pdfsync}
\usepackage{bm}
\usepackage{array}
\usepackage{amsthm,amsmath,amssymb}  
\usepackage{pgfplots}
\usepackage{epsfig,mdframed}

\usepackage{amsmath,algorithm} 
\usepackage{amssymb}  
\usepackage{color}
\usepackage{cite,shadethm,xcolor}

	\newshadetheorem{thm}{Theorem}
  \definecolor{shadethmcolor}{gray}{0.92}
	\definecolor{shaderulecolor}{HTML}{45CFFF}
	\definecolor{shaderulecolor}{gray}{0}
\newtheorem{theorem}{Theorem}
\newtheorem{lemma}{Lemma}
\newtheorem{definition}{Definition}
\newtheorem{assumption}{Assumption}
\newtheorem{proposition}{Proposition}
\usepackage{graphicx,epstopdf}

\usepackage{amsopn}

\usepackage{xspace}
\usepackage{bold-extra}
\usepackage{placeins,diagbox }

\usepackage{wrapfig} 
\usepackage{tikz,pgfplots,mathtools,float}
\usepackage{booktabs,caption}
\usepackage[flushleft]{threeparttable}
\usepackage{subfigure}
\usepackage{multirow}
\usepackage{amssymb,latexsym}  
\usepackage{extarrows}
\usepackage{dsfont} 
\usepackage{color}
\usepackage{cite}

\mdtheorem[backgroundcolor=blue!10]{set}{Setting}

\newcommand{\pmat}[1]{\begin{pmatrix} #1 \end{pmatrix}}

\usepackage[noend]{algpseudocode}
\usepackage{multirow,verbatim}
\newcommand{\Real}{\ensuremath{\mathbb{R}}}
\newcommand{\K}{{\mathbf{K}}}
\newcommand{\x}{{\mathbf{x}}}
\newcommand{\y}{{\mathbf{y}}}

\newcommand{\Kscr}{{\mathcal{K}}}
\newcommand{\Cscr}{{\mathcal{C}}}
\newcommand{\Xscr}{{\mathcal{X}}}
\newcommand{\Escr}{{\mathcal{E}}}

\newcommand{\blue}[1]{{\color{black}#1}}

\newcommand{\us}[1]{{\color{black}#1}}
\newcommand{\uss}[1]{{\color{black}#1}}

\newcommand{\ib}[1]{{\color{black}#1}}

\newcommand{\ieb}[1]{{\color{black}#1}}

\newcommand{\afj}[1]{{\color{black}#1}}
\newcommand{\cl}[1]{{\color{black}#1}}

\usepackage[colorinlistoftodos]{todonotes}

\allowdisplaybreaks
\title{Probability Maximization via Minkowski Functionals: Convex Representations and Tractable Resolution}

\author{I. E. Bardakci, A. Jalilzadeh,  C. Lagoa, and  U. V. Shanbhag\thanks{The first author reachable at \texttt{ibardakci@bartin.edu.tr}, and  the second author is reachable at \texttt{afrooz@arizona.edu} while third, and fourth authors are contactable at \texttt{cml18,udaybag@psu.edu}. The authors would like to acknowledge 
		support from NSF CMMI-1538605,  EPCN-1808266, DOE ARPA-E award DE-AR0001076,  NIH R01-HL142732, and the Gary and Sheila Bello chair funds. Preliminary efforts at studying Setting A were carried out in~\cite{badrakci18probability}. }}
\vspace{-0.7in}
\date{\today}
\begin{document}
	\maketitle
	\begin{abstract} 
		{In this paper, we consider the maximization of a probability
			$\mathbb{P}\{ \zeta \mid \zeta \in \K(\x)\}$ over a closed and convex set
			$\mathcal X$, a special case of the chance-constrained optimization problem. We define $\K(\x)$ as $\K(\x) \triangleq \{ \ieb{\zeta} \in \Kscr \mid
			c(\x,\zeta) \geq 0 \}$ where $\zeta$ is uniformly distributed on a convex and
			compact set $\Kscr$ and $c(\x,\zeta)$ is defined as either {$c(\x,\zeta)
				\triangleq  1-|\zeta^T\x|^m$, $m\geq 0$} (Setting A) or $c(\x,\zeta) \triangleq T\x -
			\zeta$ (Setting B). We show that in either setting, by leveraging recent
			findings in the context of non-Gaussian integrals of positively homogeneous
			functions,  $\mathbb{P}\{ \zeta \mid \zeta \in \K(\x)\}$ can be expressed as
			the expectation of a suitably defined \ieb{continuous} function $F(\uss{\bullet},\xi)$ with respect to an
			appropriately defined Gaussian density (or its variant), i.e. $\mathbb{E}_{\tilde p}
			[F(\x,\xi)]$. Aided by a recent observation in convex analysis, we then develop a
			convex representation of the original problem requiring the minimization of
			${g(\mathbb{E}[F(\x,\xi)])}$ over $\mathcal X$  where {$g$} is an
			appropriately defined smooth convex function. Traditional stochastic approximation schemes cannot contend with the minimization
			of ${g(\mathbb{E}[F(\bullet,\xi)])}$ over $\mathcal X$, since conditionally unbiased sampled
			gradients are unavailable. We
			then develop a regularized variance-reduced stochastic approximation ({\bf r-VRSA})
			scheme that obviates the need for such unbiasedness by combining iterative
			{regularization} with variance-reduction. Notably, ({\bf r-VRSA}) is characterized
			by both almost-sure convergence guarantees, a convergence rate of
			$\mathcal{O}(1/k^{1/2-a})$ in expected sub-optimality where $a > 0$, and 
			a sample complexity of $\mathcal{O}(1/\epsilon^{6+\delta})$ where $\delta > 0$.
			To the best of our knowledge, this may be the first such scheme for
			probability maximization problems with convergence and rate guarantees.
			Preliminary numerics on a portfolio selection problem (Setting A) and a vehicle
			routing problem (Setting B) suggest that the scheme competes well with naive
			mini-batch SA schemes as well as integer programming approximation methods.
		} 
	\end{abstract}
	\section{Introduction}
	This paper concerns the probability maximization problem ({\bf PM}), defined as 
	\begin{align} 
		\max_{\x \in \mathcal{X}} \ f(\x) \ \triangleq \  \us{\mathbb{P} \left\{\zeta \mid \zeta\in \K(\x)\right\}},
		\label{main_prob} \tag{\bf PM}
	\end{align}
	where {$f$ is a probability distribution function parametrized by a decision vector $\x$}, \us{$\mathcal{X}\subseteq \mathbb{R}^n$ denotes a closed and convex feasibility set, $\K(\x) \triangleq \{\zeta \in \mathcal{K}\mid c(\x,\zeta) \geq 0\}$, {$\mathcal{K}$ is a compact and convex set in $\mathbb{R}^n$}, $c:\mathbb{R}^n\times \mathbb{R}^d \to \mathbb{R}^m$}.  Here, $\zeta : \Omega
	\to \mathbb{R}^d$ is a $d-$dimensional random vector {with a prescribed
		distribution $\mathbb{P}$.} Problems of the form (\ref{main_prob}) fall within
	the umbrella of chance-constrained optimization problems.
	
	\subsection{Background on chance-constrained optimization} \label{sec:2.1}
	Chance-constrained optimization originates from \uss{the} 
	probabilistic  scheduling of heating oil production by Charnes,
	Cooper, and Symonds~\cite{charnes58}. A more formal description of
	chance-constrained programming as an avenue for optimization under
	uncertainty appeared in the eponymously titled paper by Charnes and
	Cooper~\cite{charnes59chance}. Such avenues have assumed relevance in
	hydro reservoir management~\cite{prekopa1978flood,ackooij14joint},
	portfolio management~\cite{pagnoncelli09saa,sun17chance}, power
	systems
	operation~\cite{bienstock2014,chanceopf1,chanceopf2,chanceuc1},
	routing~\cite{chen19routing}, structural
	failure~\cite{royset07extensions}, and inventory and supply-chain
	management~\cite{yadollahi17managing,gicquel18inventory}. \\

	\noindent \uss{{ (a) \em Analysis.} The analysis of optimization problems with probability functions has
		focused on questions of continuity,
		differentiability, and convexity. Of these, continuity and
		differentiability (and its generalized variants) are of particular
		relevance when developing algorithmic techniques. Convexity
		guarantees are important in their own right, allowing for
		certifying a stationary point as a  global maximizer. Consider a
		probability function $\psi$, defined as $\psi(\x) = \mathbb{P}[\zeta \in
		A(\x)]$ and $A: \Real^n \to \Real^m$ is a set-valued map. Under
		suitable convex-valuedness and continuity
		properties on $A$ and an appropriate measure zero requirement on $\zeta$,
		$\psi$ is continuous~\cite[Th.~2.1]{guo17convergence}. In fact, if
		$A$ is defined as $A(\x) \triangleq \{ \zeta \mid c_i(\x,\zeta)
		\leq 0, i = 1, \cdots, m\}$, then continuity of $\psi$ is implied by
		continuity of $c_i$ in both its arguments for every $i$ and a
		suitable regularity
		assumption~\cite{henrion2010optimierungsprobleme}.
		
		Differentiability of $\psi$ is a more subtle question. As eloquently
		described by van Ackooij~\cite{ackooij20review}, such results can be
		partitioned in two categories: (i) Under mild distributional
		requirements on $\zeta$ and differentiability of $c_i$ for every
		$i$, differentiability of $\psi$ can be concluded under a set of
		assumptions on $\nabla_{\zeta} c_i$, amongst others
		(cf.~\cite{uryasev92derivatives,uryasev95derivatives}); (ii)
		Alternately, by choosing the distribution, more refined statements
		are available. For instance, when the distribution of $\zeta$
		belongs to the family of eliptically symmetric distributions,
		examples being multivariate Gaussian, Student, and logistic, under
		suitable differentiability and convexity properties of $c_i$ (in the
		second argument) and additional assumptions, $\psi$ can be proven to
		be locally Lipschitz and a characterization of its Clarke
		subdifferential may be provided~\cite[Th.~1]{ackooij18subdiff}. In
		addition, if a suitable constraint qualification holds, $\psi$ can be shown to
		be differentiable and its gradient can be analytically
		characterized~\cite[Cor.~1]{ackooij18subdiff}.
		
		Finally, convexity of $\psi$ can
		be claimed under various
		conditions~\cite{prekopa72class,prekopa1973logarithmic,brascamp76extensions};
		for instance, joint quasi-convexity of $c_i$ in both arguments and
		$\alpha$-concavity of the distribution $\mathbb{P}$ implies convexity of $\psi$
		on a suitably defined set. More generally, convexity of joint chance constraints has also been studied~\cite{miller1965chance,chen2010cvar}, while more recent forays in this area have considered when the probability level is sufficiently high. Referred to as ``eventual convexity'', this avenue has been studied in the context of structured chance constraints involving copulae~\cite{van2015eventual}.\\ 
		
		\noindent { (b) \em Computation.} We now discuss the main algorithmic thrusts for resolution of chance-constrained optimization. 
		
		{(i) \underline{\em Nonlinear programming and bundle-based approaches.}} Amongst the
		earliest efforts for resolving chance-constrained optimization
		applied the penalization framework captured by the ``SUMT''
		framework, first presented by Fiacco and
		McCormick~\cite{fiacco67sumt,fiacco68nonlinear}, to the
		probabilistically constrained setting~\cite{prekopa72class}.
		Naturally, any such effort requires having deriving gradients of probability functions, as seen in the context of nonlinear probabilistic constraints with nonconvex quadratic forms~\cite{ackooij20gradient} as well as when contending with Gaussian distributions (and their variants)~\cite{ackooij14gradient,ackooij17gradient} (also see~\cite{uryasev92derivatives,uryasev94derivatives,pflug05probability}). One challenge that has been observed in early efforts is the scourge of ill-conditioning in penalizaton efforts~\cite{prekopa03probabilistic}, leading to the development of bundle-based approaches, which have proven quite powerful~\cite{ackooij14constrained,ackooij17probabilistic,ackooij21bundle}.\\

		{ (ii) \underline{\em Convex characterizations and approximations}.} 
		Convexity of (Chance-Opt) can often be claimed. For instance, ~\cite{prekopa2013}
		proved that when \us{the} distribution function of $\zeta$ is logarithmically concave (or
		log-concave) and the functions $g_1(x,\zeta)$, $g_2(x,\zeta)$, $\ldots$, $g_r(x,\zeta)$ are quasi-concave, the function $G(x)=\mathbb{P}(\zeta: g_1(x,\zeta)\geq 0, g_2(x,\zeta)\geq 0, \ldots, g_r(x,\zeta)\geq 0)$ is a log-concave function. More recently, Lagoa et al.~\cite{lagoa2005} showed that a set
		given by $\mathbb{P}\{ a^Tx \leq b\} \geq (1-\epsilon)$ is convex if $(a,b)$
		has a symmetric log-concave density and $\epsilon < 1/2$. Unfortunately, the convexity of the feasibility set does not directly allow for efficient computation, motivating a scenario-based approach. Consider (P) and its
		scenario-based approximation (P$_N$), defined next.
		$$\left\{\begin{aligned}
			\min_{\x \in \mathcal{X}} & \,  h(\x) \\
			\mbox{s.t. } & \, \mathbb{P}\left\{\zeta: c(\x,\zeta) \geq 0\right\}.
		\end{aligned} (P) 
		\right\}\xrightarrow[\mbox{\tiny Scen. approx.}]{} 
		\left\{
		\begin{aligned}
			\min_{\x \in \mathcal{X}} & \  h(\x) \\
			\mbox{s.t.} & \ c(\x,\zeta_j) \geq 0, j = 1, \cdots, N
		\end{aligned} (P_N) \right\}
		$$
		In~\cite{nemirovski2006}, the authors examined how large $N$ should be
		so that with probability $(1-\delta)$, the optimal solution of (P$_N$) is
		{\em feasible} with respect to (P) by developing conservative convex approximations. A related approach was considered by~\cite{campi11sampling}.\\ 
		
		{(iii) \em \underline{Sample-average approximation and integer programming approaches}.} Under suitable concavity assumptions on $c(\bullet,\zeta)$ and convexity requirements on $h$, (P$_N$) can be efficiently resolved. However, this is often quite conservative. Instead, the following integer program (P$_N^{\rm int}$) can serve as an approximation~\cite{luedtke08sample,ahmed17nonanticipative}.
		\begin{align}\tag{P$_N^{\rm int}$}
			\min_{\x \in \mathcal{X}, \mathbf{z} \in \{0,1\}^N} \left\{ h(\x) \mid 
			\sum_{j=1}^N z_j  \leq \gamma N, \
			\ -c(\x,\zeta_j)  \leq M_jz_j, \ j = 1, \cdots, N \right\},
		\end{align} 
		where $M_j \triangleq \max_{\x \in \mathcal{X}}
		(-c(\x,\zeta_j))$ and $\gamma \in (0,1)$. We observe that the second
		constraint ensures that only $\gamma \%$ of the $N$ scenarios are
		satisfied. In fact, when $\gamma = \epsilon$ where $\epsilon$ is the
		parameter in (P), $\hat{v}_N$ and $\hat{\x}_N$, the optimal value
		and solution of (P$_N^{\rm int}$), converge almost surely to
		$v^*(\epsilon)$ and $\mathcal{X}^*(\epsilon)$ as $N \to
		\infty$~\cite{ahmed09sample}. Naturally, if $c(\cdot,\zeta)$ is a
		nonlinear function, (P$_N^{\rm int}$) is a mixed-binary nonlinear
		program, a challenging instance of a discrete optimization problem.\\ 		
		
		{(iv) \underline{\em Smoothing-enabled
				Monte-Carlo sampling techniques.}}  Amongst the earliest
		approaches proposed by Norkin~\cite{norkin1993} utilized the characteristic function $\chi_{C}$ of a
		set $C$, defined as 
		$\chi_{C}(\zeta) = 1$ if  $\zeta \in C$ and $0$ otherwise. This allowed for  expressing $f$ as  
		$f(\x)  = \int_{\K(x)} d\mathbb{P}(\zeta) =  \int_{\mathbb{R}^d} \chi_{\K(\x)} (\zeta)d\mathbb{P}(\zeta).$
		Unfortunately, the function $\chi_C(\bullet)$ is discontinuous at the boundary of $C$, motivating the ``smoothing'' the characteristic function by using Steklov-Sobolev  smoothing~\cite{norkin1993} (also referred to as convolution-based smoothing). Specifically, $\chi_C(\bullet)$ is approximated by its smoothing $\chi_C^{\epsilon}(\bullet)$, defined as 
		$\chi^{\epsilon}_{C}(\zeta) = \int_{-\infty}^{\infty} \chi_C(\zeta+\epsilon \tau) p_s(\tau) d\tau,
		$ where $p_s(\cdot)$ is a symmetric density. The resulting approximation  $f^{\epsilon}$ is defined as 
		$f^{\epsilon}(\x)  =  \int_{\mathbb{R}^d} \int_{-\infty}^{\infty} \chi^{\epsilon}_{\K(\x)} (\zeta+\tau \epsilon) p_{s}(\tau) d\tau d\mathbb{P}(\zeta).$
		Under suitably log-concavity assumptions,
		Norkin~\cite{norkin1993} developed a stochastic approximation
		framework for maximizing the approximation $f^{\epsilon}$; However,
		there are no bounds relating the approximation and its true
		counterpart. An alternate simulation-based approach reliant on
		difference-of-convex programming~\cite{hong2011} has been recently
		proposed for chance-constrained optimization.  An alternate
		framework~\cite{campi11sampling} uses a sampling and rejection
		framework in developing estimators convergent to feasible solutions. More recent efforts have focused on obtaining stationary points of the smoothed problem~\cite{curtis18sequential,luedtke19chance}. More refined statements deriving convergence claims to Clarke stationary points have been provided in~\cite{cui20nonconvex} by conducting a variational analysis of affine chance-constrained programs.} 
	
	\subsection{Applications}\label{sec:app}
	\noindent {(1) \em \ieb{Robust portfolio selection problem.}} Portfolio selection problems
	consider the specification of  portfolio weights while  maximizing a suitable risk/reward metric while meeting risk/reward
	requirements. Much of the research in this area emerges from the seminal work by Markowitz in the 50s~\cite{markowitz52portfolio}. Consider a
	portfolio with $n$ risky assets, whose random returns are denoted by a random
	variable $\ib{\zeta=[\zeta_1,\zeta_2,\ldots,\zeta_n]^T}$ with mean returns $\mu
	= [\mu_1,\mu_2, \cdots, \mu_n]^T$ and covariance $ \Sigma $. Let the proportion of the portfolio to be invested in each asset be represented by $\x=[x_1,x_2,\ldots,x_n]^T$. It is
	assumed that no assets will be shorted and, hence, without loss of generality,
	the set of all possible portfolio allocations is given by $\mathcal{X}
	\triangleq \left\{\x: {\bf 1}^T \x =1 \text{ and } \x \geq 0 \right\}$.  \ieb{We consider the robust portfolio selection problem (RPS) where the distribution of asset returns is not known but some of its properties available.} Given
	a threshold $\alpha$ and an allocation $\x$, \ieb{the distributionally robust risk associated with portfolio weights $ \x $ is defined as $f_{\alpha}(\x) \triangleq \sup_{\upsilon \in \mathcal{H}}\mathbb{P}_{\upsilon}
		\left\{\zeta: \zeta^T\x\leq -\alpha\right\}$ where $  \mathbb{P}_{\upsilon}$ denotes probability computed using the probability density function $ \upsilon $ belonging to admissible class of density functions $ \mathcal{H}.$} 
	\cl{This class of distributions contains all distributions whose density is radially decreasing and have level sets of the form $\{\delta: \delta^T \Sigma^{-1} \delta = r\}$ for some $r$.}
	Consider the following \ieb{distributionally robust} portfolio selection problem with $\gamma > 0$.
	\begin{align} \tag{RPS}
		\min_{\x\in \mathcal{X}} \left\{\sup_{\upsilon \in \mathcal{H}}\mathbb{P}_{\upsilon}
		\left\{\zeta: \zeta^T\x\leq -\alpha\right\} \mid  \mu^T \x = \gamma \right\}.
	\end{align}
	\ieb{
		In prior work~\cite{bardakci2019distributionally}, it was shown via the following lemma that the supremum in (RPS) is achieved when $ \upsilon $ is a uniform density over an ellipsoid.
		\begin{lemma} \label{lem: worst_dist}	
			Let the random vector $ \zeta $ is of the form $ \zeta = 
			\boldsymbol{\mu} + \Delta  $ where the distribution~$h$ for~${\Delta}$ is taken to be unknown but assumed to belong to the class  $\mathcal{H}$ and
			$  	\boldsymbol{\mu} \in \{\zeta : \zeta^T \x \geq -\alpha \}.$	Then \uss{$ 	\sup_{\upsilon \in \mathcal{H}}  \text{Prob}_{\upsilon}\{\zeta : \zeta^T \x \leq -\alpha \} $} is achieved when $ \Delta$ has a uniform distribution over the set
			\begin{equation} \label{set:R}
				\mathcal{R} = \{ \Delta \in \mathbb{R}^n:  \Delta^T \Sigma^{-1}  
				\Delta \leq r_{\max}\},
			\end{equation}
			where $ r_{\max} $ represents the uniform bound on the support of the  
			distribution of $\Delta$. 
		\end{lemma}
		Consequently, the robust portfolio selection problem (RPS) reduces to the probability minimization problem. 
		\[
		\min_{\x\in \mathcal{X}} \left\{\mathbb{P}_{\upsilon}
		\left\{\zeta: \zeta^T\x\leq -\alpha\right\} \mid  \mu^T \x = \gamma \right\},
		\]
		where $\mathbb{P}_{\upsilon}$ denotes the uniform distribution over an ellipsoid. This motivates the consideration of such a portfolio selection problem in the numerics. }

	\noindent {(2) \em Set covering problems.} \us{Consider a set covering problem~\cite{shapiro2009} (\uss{closely related to a vehicle routing problem}) 
		\begin{align} \label{ex1}
			\max \  \left\{f(\x) \mid \x \in \mathcal{X}\right\}, \mbox{ where }  
			f(\x) \ \triangleq \ \mathbb{P}\{\zeta \in \Kscr \mid T \x \geq \zeta\},
		\end{align}
		where $\mathcal{X} \triangleq \left\{\x \mid c^T\x \leq \ieb{\beta}, \x \geq 0\right\}$, {and $T\in \mathbb R^{d\times n}.$}} The incidence matrix $ T $ represents a network \afj{with $d$ arcs} and $n$ routes {with $ij$th component denoted by $t_{ij}$ where  
		$t_{ij} \triangleq 
		1$ if route $ j $ contains arc $ i $ and is zero otherwise}. Furthermore, $\ib{\zeta_i}$ denotes the random demand on arc $i$  where $\ib{\zeta \triangleq (\zeta_1, \cdots, \zeta_d)^T }\in \Kscr$, $c_j$ represents the non-negative cost of operating route $j$, and $\ieb{\beta}$ represents a given cost threshold.

	\subsection{Gaps, contributions, and outline} {\em {{Gaps}.}} The
	optimization of distribution functions remains a 
	challenging problem. To the best of
	our knowledge, {\em there exist no efficient schemes equipped with 
		asymptotic convergence or  rate guarantees for probability maximization problems or their generalizations, i.e. chance-constrained problems}. This is both a testament to the difficulty of such a problem as well
	as a motivation for the present work which intends provides precisely
	such schemes with suitable convergence and rate guarantees for a
	subclass of problems. The key contributions of this work are as
	follows.

	\noindent {\em (I) Representation of $(\ref{main_prob})$ as a convex program.}
	By leveraging recent findings on non-Gaussian integrals of positively
	homogeneous functions (PHFs) \cite{lassere2015,morozov2009}, we consider regimes where $c(\x,\zeta) = 1-
	\ib{|\zeta^T}\x|^2$ (Setting A, Section~\ref{sec:2.2}) {or $c(\x,\zeta) = T\x -\zeta$ (Setting
		B, Section~\ref{sec:2.3})} where $\Kscr$ is a compact and convex set, symmetric about the origin. In both settings, we show that $\mathbb{P}\{\zeta
	\mid \zeta \in \K(\x)\}$ can be expressed as the expectation of a suitably
	defined \uss{Clarke regular} integrand $F(\x,\xi)$, i.e. $ \mathbb{E}_{\tilde p}
	[F(\x,{\xi})]$ where \blue{$\tilde p(\ib{\xi})$} is either a suitably defined Gaussian density or its variant.
	We then
	proceed to show that the original problem is equivalent to the minimization of
	a convex function $g(\mathbb{E}_{\tilde p}[F(\bullet,\xi)])$ over a closed and convex set $\mathcal X$ where $g$ is a suitably defined
	convex and smooth function. We refine these relationships when $\Kscr$ is
	either an $\ell_p$ ball or an ellipsoid centered at the origin (Setting A) or
	loses symmetry (Setting B).
	
	
	\noindent {\em (II) Regularized variance-reduced stochastic approximation {\bf (r-VRSA)} scheme.} The
	resulting convex program is an instance of a compositional stochastic
	optimization problem where unbiased first-order oracles are unavailable. In
	Section~\ref{sec:3}, we present a regularized variance-reduced SA
	scheme that combines variance-reduction (to accommodate bias) with
	regularization. The resulting scheme is characterized by a
	rate of convergence of $\mathcal{O}(1/{k^{(1/2-a)}})$ and
	$\mathcal{O}(1/K^{1/2})$ for diminishing and constant steplengths (the latter requiring the specification of a simulation length $K$) while the
	sample-complexity to achieve $\epsilon-$optimality, i.e.
	$\mathbb{E}[h(\x^*)-h(\x_K)] \leq \epsilon$ is $\mathcal{O}(1/\epsilon^{6+\delta})$ and $\mathcal{O}(1/\epsilon^6)$, respectively.
	While the rate of convergence matches the optimal rate for subgradient methods
	for convex programs, the sample-complexity is worse than the canonical
	$\mathcal{O}(1/\epsilon^2)$. The latter is unsurprising since we do not have
	access to unbiased oracles. {It is worth emphasizing that this appears to be
		one of the first schemes with asymptotic convergence and rate guarantees for a
		class of non-trivial probability maximization problems.} Apart from this, we
	provide some background in Section~\ref{sec:2.1} and numerical experiments are
	discussed in Section~\ref{sec:4}.
	
	{\bf Notation.}  We conclude this section with a review of {\em notation.} The sets of real numbers, {non-negative real numbers}, non-negative integers, and positive integers
	are denoted by $\mathbb{R}$, {$\Real_+$}, $\mathbb{N}$, and $\mathbb{Z}$,
	respectively. The Euclidean norm of column vectors $\mathbf{x} \in
	\mathbb{R}^n$ is denoted by $\|\mathbf{x}\|$, while the spectral norm of $\mathbf{A} \in \mathbb{R}^{m\times n}$ is given by $\| \mathbf{A} \| = \text{max}\{ \|\mathbf{Ax} \| \colon \|\mathbf{x}\| \leq 1 \}$. The $n$-by-$n$ identity matrix is written as $\mathbf{I}_n$, and the $m$-by-$n$ zero matrix as $\mathbf{0}_{m \times n}$. The projection onto the set $X$ is denoted by $\Pi_{X}$, that is, $\Pi_{X}(y) = \text{argmin}_{x \in X} \| x - y \| $. {Finally, unless mentioned otherwise, any missing proofs are provided in the appendix.} 
	
	\section{\us{An expectation-valued convex framework}}\label{sec:2}
	{Throughout this paper, we consider (\ref{main_prob}) where
		\begin{align}\label{def-K}
			\K(\x) \triangleq \left\{ \zeta \in \mathcal{K} \colon c(\x,\zeta) \geq 0 \right\}. 
		\end{align}
		We consider two sets of regimes based on the choice of $\mathcal{K}$ and $c(\x,\zeta)$. In addition, we impose an assumption on $\Kscr$ and a distributional assumption on $\zeta$.
		\begin{assumption}[Assumptions on $\Kscr$ and $\mathcal X$] \label{dist-ass} The random variable $\zeta$ is uniformly distributed on the set $\Kscr$ where $\Kscr$ is a compact and convex set in $\Real^n$ and is symmetric about the origin. The set $\mathcal X$ is closed, convex, and bounded. 
		\end{assumption}  
		\begin{set*}[A]	
			The constraint $c(\x,\zeta)$ in \eqref{def-K} is defined as
			\begin{align}  \label{type1}
				c(\x,\zeta) \triangleq \us{1-|\zeta^{\intercal}\x|^m}, 
			\end{align}	
			where $m  \in \Real_+$ and $\zeta \in \Real^n$ is a random variable.    \qed 
		\end{set*}
		\begin{set*}[B]  The constraint $c(\x,\zeta)$ in \eqref{def-K} is defined as 
			\begin{align}  \label{type2}
				c(\x,\zeta) \triangleq T\x -\zeta,
			\end{align}
			where $T \in \Real^{d \times n}$ and $\zeta \in \Real^d$ is a random variable. 
			\qed
		\end{set*}  
		We observe that Setting A can capture problems such as the portfolio optimization problem described in Section~\ref{sec:app}. \us{Without loss of generality, $m = 2$ in this paper.}   Further, in such problems, $f({\bf 0}) = 1$ and 
		$\lim_{\|\x\| \to \infty} f(\x) = 0, $ where $f(\x) \triangleq \mathbb{P}(\{\zeta \mid \zeta \in \K(\x)\}.$ Setting B assumes relevance when considering chance constraints defined using polyhedral constraints with uncertain right-hand sides. An instance of such a problem is the \uss{set covering} problem described in Section~\ref{sec:app}.  We qualify problems in Settings A and B as (\ref{main_prob}$_A$) and (\ref{main_prob}$_B$), respectively.  Before proceeding, we recall two definitions of relevance. 
		\begin{definition} [Log-concavity, positive homogeneity]
			A function $f : \mathbb{R}^d \rightarrow [0, \infty)$ is 
			log-concave if for any $x, y \in \mathbb{R}^n $ and $\lambda \in [0,1]$, 
			$f((1-\lambda)x + \lambda y) \geq [f(x)]^{1-\lambda} [f(y)]^{\lambda}.$
			A continuous function $ f: \mathbb{R}^{\cl{d}} \rightarrow \mathbb{R} $ is called positively homogeneous function of degree $ p \in \mathbb{R} $ if $ f(\alpha x) = \alpha^p f(x)$ for all $ \alpha>0 $ and all $ x \in \mathbb{R}^{\ieb{d}} $.
			
		\end{definition}
		\begin{definition} [Minkowski Functional]
			Let the set $\Kscr \subset \mathbb{R}^n$. Then, the Minkowski functional associated with the set $\Kscr$, denoted by $\|\cdot\|_\Kscr$, is defined as 
			$\|\zeta\|_{\Kscr} \triangleq \inf \{t>0 : \zeta/t \in \Kscr\}$
			for all $\zeta \in \mathbb{R}^n $. 
		\end{definition}
		Note that $\|\cdot\|_{\Kscr}$ defines a norm when $\Kscr$ is compact,
		convex and symmetric.  For instance, if $\Kscr$ is the unit ball in
		$\mathbb{R}^n$, then the Minkowski functional reduces to $\|\cdot\|_2$ in
		$\mathbb{R}^n$ i.e. $\|\zeta\|_{\Kscr} = \|\zeta\|_2.$ In the remainder of
		this section, after providing some background in Section~\ref{sec:2.1}, we
		proceed to show that the function $f$, \us{defined in \eqref{main_prob}},
		is equivalent to the expectation of a nonsmooth integrand. In fact, \uss{we prove that this integrand is Clarke regular and}  the reciprocal of $f$ is a convex function for setting A
		(Section~\ref{sec:2.2}) while \us{the negative log-transformation of $f$ is
			a convex function in setting B (Section~\ref{sec:2.3}).} 
		
		\subsection{Expectation-valued convex representations for Setting A.}\label{sec:2.2} 
		We begin by recalling that  $\mathbb{P}\{\zeta \mid \zeta \in \K(\x) \}$ can be rewritten as 
		\begin{align}
			f(\x) =  \mathbb{P} \{\zeta \mid \zeta \in \K(\x)\} =  \frac{1}{\text{Vol}(\mathcal{K})} \int_{\K(\x)}
			{1} \ d\zeta, \label{uni}
		\end{align}
		where \us{the last equality follows from Assumption~\ref{dist-ass}} and $\text{Vol}(\mathcal{K})$ denotes the volume of the set $\mathcal{K}$.  We now show that \eqref{uni} can be expressed as an expectation with respect to
		prescribed probability measure.   
		Recall that a function  $\uss{r}$ is a PHF of degree $m$ where $\uss{r}(\zeta) = \max\{\uss{r}_1(\zeta),\hdots, \uss{r}_\ell(\zeta)\}$ if $\uss{r}_1, \hdots, \uss{r}_\ell$ are PHFs of degree $m$. Lemma~\ref{las-1} provides conditions under which the integral of a PHF over a suitable set is equal  to another integral which is an expectation over a suitably defined measure.
		\begin{lemma}\label{las-1}{\em\cite[Cor.~2.3]{lassere2015}}
			{Let $h$ be a positively homogenous function of degree $p$ and}
			let $\ieb{r}_1,\hdots,\ieb{r}_\ell$ be positively homogeneous functions (PHFs)
			of degree $0 \neq \us{t} \in \mathbb{R}$. Let $\uss{\Lambda}$ be a bounded set defined as $\uss{\Lambda} \ \triangleq \
			\{\zeta: \ieb{r}_k(\zeta) \leq 1, k=1,\hdots, \ell\}.$ 
			If $ \int_{\mathbb{R}^n} \us{\left|h(\xi)\right| } e^{-\us{\max\{\ieb{r}_1(\xi), \hdots, \ieb{r}_{\ell}(\xi)\}}} \ d\xi < \infty,$ then the following holds.
			\begin{align*}
				\int_{\uss{\Lambda}} \us{h(\zeta)} \ d\zeta = \frac{1}{\Gamma(1+\us{(n+p)}/\us{t})} \int_{\mathbb{R}^n} \us{h(\xi)} e^{-\us{\max\{\ieb{r}_1(\xi), \hdots, \ieb{r}_{\ell}(\xi)\}}} \ d\xi.
			\end{align*} 
		\end{lemma}
		{We now show that $f$ is given by the expectation with
			respect to  $\tilde p(\xi)$, the density function of a suitably defined random variable dependent on the choice of $\Kscr$.}  
		\begin{theorem}[Representation of $f$ as expectation for general $\Kscr$] \label{theo-rep} Consider  \eqref{main_prob}. Suppose Assumption~\ref{dist-ass} holds,  where $c$ is defined as (\ref{type1}) and $\K$ is defined as \eqref{def-K}. Then the following equality holds.
			\begin{align}
				\label{type1-exp} & \ \mathbb{P} \left\{\zeta: \zeta \in \K(\x)\right\} = \mathbb{E}_{\tilde p(\xi)} [F(\x,\xi)] = \int_{\mathbb{R}^n} F(\x,\xi) \tilde p(\xi) d\xi,\mbox{ where } \\ 
				\label{def-F-A} F(\x,\xi) & \triangleq  \mathcal{C}_{\Kscr}(2\pi)^{n/2} 
				e^{-\max\{|\xi^{\intercal} \x |^2, {\|\xi\|^2_{\mathcal{K}}}\}+\frac{\|\xi\|^2_{\Kscr}}{2}}, 
				\mathcal{C}  \triangleq \tfrac{1}{\mathrm{Vol}(\Kscr)} \afj{\tfrac{1}{\Gamma(1+n/\us{2})}},\\
				\label{def-p-A}\tilde p(\xi) &  \triangleq 
				\tfrac{1}{ (2\pi)^{n/2} D_{\Kscr} } e^{\frac{-\|\xi\|^2_{\Kscr}}{2}},
			\end{align} 
			\us{$D_{\Kscr}$ is a positive scalar such that $\int_{\Real^n}  
				\tfrac{1}{ (2\pi)^{n/2} D_{\Kscr} } e^{\frac{-\|\xi\|^2_{\Kscr}}{2}} = 1$, and $\mathcal{C}_{\Kscr} \triangleq \ieb{\mathcal C} D_{\Kscr}$.}
		\end{theorem}
		\begin{proof}
			{We begin by noting that $\K(\x)$ can be expressed as:
				$$ \K(\x) = \{\zeta: \zeta \in \mathcal{K}\} \cap \{\zeta: |\zeta^{\intercal} \x | \leq 1\} $$}  
			Since the set $\mathcal{K}$ is
			convex, compact, and symmetric, the Minkowski functional of
			$\mathcal{K}$ defines a norm, and hence, it is a PHF. Moreover, by
			the definition of the Minkowski functional we have
			$\zeta \in \mathcal{K} \Leftrightarrow  \|\zeta\|_{\Kscr} \leq 1.$
			By using this definition, we may rewrite $\K(x)$ as follows.
			\begin{align*}
				\us{\K(\x)} =  \{\zeta: |\zeta^{\intercal} \x |^{\us{2}} \leq 1\} \cap \{\zeta
				\colon \|\zeta\|^{\us{2}}_\mathcal{K} \leq 1\} =  \left\{\zeta: \text{max}\{|\zeta^{\intercal} \x |^\ieb{2}, \|\zeta\|^{\us{2}}_\Kscr\} \leq 1
				\right\}.
			\end{align*}
			\us{Since $|{\zeta}^\intercal \x |^{\us{2}}$ and $ \|\zeta\|^{\us{2}}_{\mathcal{K}}$ are both PHFs of
				degree $\us{2}$, then \afj{$g(\x,\uss{\bullet})$} is also a PHF of degree $\us{2}$ where \afj{$g(\x,\zeta)$} is defined as 
				$\afj{g(\x,\zeta)} \triangleq \text{max}\{|\zeta^\intercal \x |^\us{2},
				{\|\zeta\|^{\us{2}}_\mathcal{K}}\}$.} By selecting $h(\zeta) = 1$ and $\ieb{\Lambda} = \K(\x)$, we may invoke Lemma~\ref{las-1}, leading to the following equality. 
			\begin{align} 
				f(\x)= \frac{1}{\text{Vol}(\mathcal{K})} \int_{\us{\K(x)}} 1 \  d\zeta = \frac{1}{\text{Vol}(\mathcal{K})} \frac{1}{\Gamma(1+n/\us{2})} \int_{\mathbb{R}^n} e^{\afj{-g(\x,\xi)}} \ d\xi, \label{Las}
			\end{align}
			whenever $\int_{\mathbb{R}^n} e^{-g(\ieb{\x},\xi)} \ d\xi$ is finite. In fact, the expression (\ref{Las}) can be written as
			\begin{align}
				f(\x)  \ & = \mathcal{C}  \int_{\mathbb{R}^n} \left( (2\pi)^{n/2} e^{-\max\{|\xi^{\intercal} \x |^{\us{2}}, {\|\xi\|^{\us{2}}_{\mathcal{K}}}\}+\frac{\|\xi\|^2_{\Kscr}}{2}}\right)\left( \tfrac{1}{(2\pi)^{n/2}}\afj{e^{\frac{-\|\xi\|^2_{\Kscr}}{2}}} 
				\right)  \ d\xi \nonumber \\ 
				\ & =  \int_{\mathbb{R}^n} \underbrace{\left(\mathcal{C}_{\Kscr} (2\pi)^{n/2} e^{-\max\{|\xi^{\intercal} \x |^{\us{2}}, {\|\xi\|^{\us{2}}_{\mathcal{K}}}\}+\frac{\|\xi\|^2_{\Kscr}}{2}}\right)}_{\triangleq F(\x,\xi)}\underbrace{\left( \tfrac{1}{D_{\Kscr}(2\pi)^{n/2}} \afj{e^{\frac{-\|\xi\|^2_{\Kscr}}{2}}} 
					\right)}_{\triangleq \tilde p(\xi)}  \ d\xi \nonumber \\ 
				\notag & =   \int_{\mathbb{R}^n} F(\x,\xi) \ \tilde p(\xi) \ d\xi  = \mathcal{C} \ \mathbb{E}_{\tilde p}[F(\x,\xi)], 
				\mbox{ where } \mathcal{C}  \triangleq \tfrac{1}{\text{Vol}(\Kscr)} \afj{\tfrac{1}{\Gamma(1+n/\us{2})}}, 
			\end{align} 
			$\tilde p(\xi)$  denotes the density, 
			$D_{\Kscr}$ is such that $\int_{\Real^n} \tilde{p}(\xi) d\xi = 1$, and $\mathcal{C}_{\Kscr} \triangleq \mathcal{C} D_{\Kscr}$. 
		\end{proof}
		Next, we examine the convexity properties of a related problem, given by (\ref{hx}). 
		\begin{align}
			\min_{\x \in \mathcal{X}} \ h(\x)  \  \triangleq \   \frac{1}{f(\x)}. \label{hx}
		\end{align}
		{Crucial to this claim is the leveraging of a result provided in~\cite{bobkov2010} which allows for claiming that the reciprocal of \afj{$\mathbb{P}\{\zeta \in \K(\x)\}$} is a convex function in $\x$ when $\zeta$ satisfies a suitable requirement.} 
		\begin{theorem}[Transformation of (\ref{main_prob}$_A$) to convex program]  \label{theo-convexity-settingA} Suppose  the function $f$, $\K(\x)$, and $\mathcal{X}$ are as defined in \eqref{main_prob}. Suppose $f(\x) \in [\epsilon,1]$ for $\x \in \mathcal X$ \afj{and $\epsilon>0$} and $h$ is defined such that $h(\x) = 1/f(\x)$. Then the following hold. 
			\begin{enumerate}
				\item[(a)] $h$ is convex in \us{$\x$ over $\mathcal X$} where $h(\x) \triangleq 1/f(\x)$. 
				\item[(b)] A global maximizer of \eqref{main_prob} is a global minimizer of \eqref{hx}.
			\end{enumerate}
		\end{theorem} 
		
		Before proceeding, we provide a lemma for computing the maximal value of $u^ce^{-u}$ where $c$ is a positive integer. 
		\begin{lemma}\label{unimodal}
			Consider the function $u^ce^{-u}$ defined on $u \in \mathbb{R}_+$ where $c \ge1 $ and $c \in \mathbb{Z}_+$. Then we have that 
			$ \displaystyle \max_{u \geq 0} u^ce^{-u} = \tfrac{c^c}{e^c}$ and $\mbox{\em arg}\displaystyle \max_{u \geq 0} u^c e^{-u} = c.$ 
		\end{lemma}
		
		\uss{Note that since $F(\bullet,\xi)$ is not necessarily convex, we cannot
			employ subdifferentials of $F$. Instead, we  begin by recalling some key elements of Clarke's
			nonsmooth calculus and start by providing  the definition of the
			Clarke generalized gradient of a function $h$ by leveraging its directional
			derivatives. 
			
			\begin{definition}[{\bf Directional derivatives and Clarke generalized gradient~\cite{clarke98}}] \em 
				The directional derivative of $h$ at $\x$ in a direction $v$ is defined as 
				\begin{align}
					h^{\circ}(\x,v) \triangleq  \limsup_{\y \to \x, t \downarrow 0} \left(\frac{h(\y+tv)-h(\y)}{t}\right).
				\end{align}
				The Clarke generalized gradient at $\x$ can then be defined as 
				\begin{align}
					\partial h(\x) \triangleq  \left\{ \xi \in \Real^n \mid h^{\circ}(\x,v) \geq  \xi^T v, \quad \forall v \in \Real^n\right\}.
				\end{align}
				In other words, $h^{\circ}(\x,v) = \displaystyle \sup_{g \in \partial h(\x)}  g^Tv.$
			\end{definition}
			If $h$ is $C^1$ at $\x$, the Clarke generalized gradient reduces to the standard gradient, i.e. \cl{$\partial h(\x)$ is a singleton at $\nabla_{\x} h(\x).$}  We now  review some properties of $\partial h(\x)$. In particular, if
			$h$ is locally Lipschitz on an open set $\Cscr$ containing $\Xscr$, then $h$ is
			differentiable almost everywhere on $\Cscr$ by Rademacher's
			theorem~\cite{clarke98}. Suppose $\Cscr_h$ denotes the set of points where $h$
			is not differentiable. We now provide some properties of the Clarke generalized
			gradient.   \begin{proposition}[{\bf Properties of Clarke generalized
					gradients~\cite{clarke98}}] \em
				Suppose $h$ is $L_0$-Lipschitz continuous on $\Real^n$. Then the following hold.
				\begin{enumerate}
					\item[(i)] $\partial h(\x)$ is a nonempty, convex, and compact  set and $\|\uss{u} \| \leq {L_0}$ for any $\uss{u} \in \partial h(\x)$. 
					\item[(ii)] $h$ is differentiable almost everywhere. 
					\item[(iii)] $\partial h$ is an upper semicontinuous map defined as 
					$$\hspace{-.5in}\partial h(\x) \triangleq \mbox{conv}\left\{\uss{u} \mid \uss{u} = \lim_{k \to \infty} \nabla_{\x} h(\x_k), \Cscr_h \not \owns \x_k \to \x\right\}.  
					$$
				\end{enumerate}
		\end{proposition}}
		
		\uss{To employ the Clarke generalized gradient, we require that the function be at least locally Lipschitz. We proceed to prove that $F(\bullet,\xi)$ satisfies this requirement. We further show that this result paves the way for showing that we can interchange the Clarke subdifferential and the expectation operator.}
		\begin{lemma}\label{lemma_Clarke_F}
			Consider the function $F(\bullet,\xi)$ defined as
			\begin{align*}
				F(\x,\xi) = \begin{cases}  \left( \ieb{\mathcal{C}_{\Kscr}}(2\pi)^{n/2} e^{-|\xi^{\intercal} \x |^{\us{2}}+\frac{\us{\|\xi\|_{\Kscr}^2}}{2}}\right)  & \xi \in \Xi_1(\x) \triangleq \left\{\xi \mid |\xi^{\intercal} \x |^{\us{2}} > \|\xi\|^{2}_{\mathcal{K}}\right\} \\ 
					\left( \ieb{\mathcal{C}_{\Kscr}}(2\pi)^{n/2} e^{-\max\{|\xi^{\intercal} \x |^{\us{2}}, {\|\xi\|^{\us{2}}_{\mathcal{K}}}\}+\frac{\us{\|\xi\|_{\Kscr}^2}}{2}}\right)  & \xi \in \Xi_0(\x) \triangleq \left\{\xi \mid |\xi^{\intercal} \x |^{\us{2}} = \|\xi\|^{2}_{\mathcal{K}}\right\}\\
					\left( \ieb{\mathcal{C}_{\Kscr}} (2\pi)^{n/2} e^{-\frac{\us{\|\xi\|_{\Kscr}^2}}{2}}\right).  & \xi \in \Xi_2(\x) \triangleq \left\{\xi \mid |\xi^{\intercal} \x |^{\us{2}} < \|\xi\|^{2}_{\mathcal{K}}\right\}
				\end{cases}
			\end{align*}
			Then the following hold.
			
			\noindent (a)  $F(\bullet,\xi)$ is locally Lipschitz  for every $\xi$. 
			
			\noindent (b) $F(\bullet,\xi)$ is a Clarke regular function for almost every $\xi \in \Real^n$. 
			
			\noindent (c) For any $\x \in \Real^n$, $ \partial \mathbb{E}[F(\x,\xi)] = \mathbb{E}[\partial F(\x,\xi)].$    
		\end{lemma}
		\begin{proof} (a) This follows by observing that  $F(\bullet,\xi)$ is C$^1$ when $\xi \in \Xi_1(\x) \cup \Xi_2(\x)$ and  piecewise C$^1$ if $\xi \in \Xi_0(\x)$.  Therefore $F(\bullet,\xi)$ is locally Lipschitz for every $\xi \in \Xi$~\cite[Cor.~4.1.1.]{scholtes2012introduction}.
			
			\noindent (b) Since $\Xi_0(\x)$ is a lower-dimensional set in $\Real^n$, we have that $\Xi_1(\x) \cup \Xi_2(\x) = \Real^n$. Therefore for almost every $\xi \in \Real^n$, we have that $F(\bullet,\xi)$ is C$^1$. Consequently, $F(\bullet,\xi)$ is a Clarke regular function for almost every $\xi$.
			
			\noindent (c) Since $F(\bullet,\xi)$ is Clarke regular  for almost every $\xi \in \Real^n$, by ~\cite[Theorem 3.4.]{burke20subdifferential}, we have that 
			$ \partial \mathbb{E}[F(\x,\xi)] = \mathbb{E}[\partial F(\x,\xi)].$    
		\end{proof}
	} 
	
	Computational schemes, particularly via stochastic approximation, rely on
	boundedness of $F(\bullet,\xi)$ and \uss{$G(\x,\xi)$ where  we denote an element of $\partial F(\x,\xi)$ by $G(\x,\xi)$, i.e. $G(\x,\xi) \in \partial F(\x,\xi)$}.
	
	\begin{proposition}[Properties of $F(\x,\xi)$ and $G(\x,\xi)$ under general $\Kscr$]\label{bound sub}
		Consider the function $f$ in \eqref{main_prob} and suppose Assumption~\ref{dist-ass} holds. Suppose  $c(\x,\zeta)$, $F(\x,\xi)$, and $\tilde{p}(\xi)$ are  defined as (\ref{type1}), \eqref{def-F-A}, and \eqref{def-p-A}, 
		respectively and $G(\x,\xi) \in \partial F(\x,\xi)$ for any $\x \in \mathcal{X}$. Then the following hold.
		
		\noindent (a) For any $\x \in \mathcal{X}$, $|F(\x,\xi)|^2 \leq \mathcal{C}^2_{\Kscr} (2\pi)^n$ for every $\xi \in \Real^n$.
		
		\noindent (b) For any $\x \in \mathcal X$, 
		$\mathbb{E}[\|G(\x,\xi)\|^2] \leq \frac{\mathcal{C}^2_{\Kscr}(2\pi)^n}{e} \mathbb{E}_{\tilde p}[\|\xi\|^2].$ 
		
	\end{proposition}
	
	\us{We may specialize this \ieb{representation and the} bounds to regimes where $\Kscr$ is an $\ell_p$-ball in $\Real^n$ via the following Lemma. In addition, we recall that the density of a multivariate Gaussian with independent components, each with mean zero and variance $\sigma^2$, has a density given by $\tilde p(\xi)$ defined as 
		$$ \tilde p(\xi) \triangleq \tfrac{e^{-\frac{\|\xi\|^2_2}{2\ieb{\sigma^2}}}}{(2\pi \sigma)^{n/2}}, \mbox{ where } \int_{\Real^n} \tilde{p}(\xi) d\xi = 1.$$  
		Consequently, $D_{\Kscr} = 1$ and $\mathcal{C}_{\Kscr} = C$. We rely on the following standard lemma.}
	\begin{lemma}\label{cor: norm_ineq}
		Let  the $\ell_p$-norm \us{of a vector $x \in \mathbb{R}^n$} be defined as $ \|x\|_{p} \triangleq \left(\sum_{i=1}^{n} \mid x_i\mid^{p} \right)^{1/p}$. For any $ 1\leq a<b $, {there exists a scalar $\beta \triangleq n^{(1/a-1/b)}$ such that for every $x\in \Real^n$, 
			$\|x\|_{b} \leq \|x\|_{a} \leq \beta \|x\|_{b}.$}
	\end{lemma}

	\begin{proposition}[Representation and boundedness of $F(\x,\xi)$ and $G(\x,\xi)$ when $\Kscr$ is an $\ell_p$ ball] \label{bound-lp-ball-sett_A}
		Consider the function $f$ in \eqref{main_prob}. Suppose Assumption~\ref{dist-ass} holds and $c(\x,\xi)$ is defined as (\ref{type1}). Suppose $\Kscr$ is an $\ell_p$-ball in $\Real^n$ where $p \geq 1$.  Then the following hold.
		
		\noindent (a) $f(\x) = \mathbb{E}_{\tilde p}[F(\x,\xi)]$, where
		\begin{align*}
			F(\x,\xi) &  \triangleq \left( \mathcal{C} (2\pi \sigma^2)^{n/2} e^{-\max\{|\xi^{\intercal} \x |^{\us{2}}, {\|\xi\|^{\us{2}}_{p}}\}+\frac{\us{\|\xi\|_{2}^2}}{2\sigma^2}}\right),\ieb{\mathcal{C}\triangleq \tfrac{1}{\mathrm{Vol}(\Kscr)} \tfrac{1}{\Gamma(1+n/2)}}, \mbox{ and},\\
			\tilde p(\xi) &\triangleq \tfrac{1}{(2\pi \sigma^2)^{n/2}} e^{-\frac{\|\xi\|_{2}^2}{2 \sigma^2}},  \
			\mbox{ where } \sigma^2  \triangleq \begin{cases}
				n^{1/2-1/p}, & p \geq 2 \\
				1. & 1 \leq p < 2
			\end{cases}
		\end{align*}
		
		\noindent (b) For any $\x \in \mathcal{X}$, $|F(\x,\xi)|^2 \leq \ieb{\mathcal{C}}^2 (2\pi \sigma^2)^n$ for every $\xi \in \Real^n$.

		\noindent (c) For every $\x \in \mathcal X$,  $\mathbb{E}[\|G(\x,\xi)\|^2] \leq  e^{-1} \ieb{\mathcal{C}}^2(2\pi\sigma^2)^{n}\mathbb{E}_{\tilde p}[\|\xi\|^2],$
		where $G(\x,\xi) \in \partial F(\x,\xi)$ for any $\x \in \mathcal{X}$ \uss{and $\xi \in \Real^n$}.
	\end{proposition}
	Next, we consider the regime where $\Kscr \triangleq \Kscr_{\Escr}$ is an ellipsoid in $\Real^n$,  defined as 
	\begin{align}\label{ellipsoid} \Kscr_{\Escr} \triangleq \left\{ \zeta \in \Real^n \mid \zeta^{\intercal} U^{\intercal} \ib{\Sigma^{-1}} U \zeta \leq 1\right\}, \end{align}
	where $U \in \Real^{n \times n}$ is an orthogonal matrix whose columns represent unit vectors along the principal axes of the ellipsoid and $\ib{\Sigma^{-1}}$ is a positive diagonal matrix with the $i$th diagonal element denoted by $\ib{1/ \sigma^2_{ii}}$. By defining $\ib{\eta = \Sigma^{-1/2} U\zeta}$ or $U^{\intercal} \ib{\Sigma^{1/2} \eta = \zeta}$, we may observe that $\Kscr_{2} = \ib{\{\eta: \|\eta\|^2_2 \leq 1\}}$, i.e. $\Kscr_{\Escr}$ can be transformed to $\Kscr_2.$ We now prove that there is an equivalence between (\ref{main_prob}$\uss{^{\Escr}_A})$ and its \uss{transformed} counterpart (\ref{main_prob}$\uss{^{2}_A}$). 
	\begin{lemma}[Equivalence between (\ref{main_prob}$^{\Escr}_A$) and (\ref{main_prob}$^{2}_{A}$)] \label{equiv-ellipsoid}  
		Consider the function $f$ in \eqref{main_prob}. Suppose Assumption~\ref{dist-ass} holds and $c$ is defined as (\ref{type1}). Suppose $\Kscr_{\Escr}$ is an ellipsoid in $\Real^n$, \afj{ defined in \eqref{ellipsoid}}.
		Then \uss{$\x$} is a solution of $($\ref{main_prob}$^{\Escr}_{A})$ 
		\begin{align} \tag{\ref{main_prob}$^{\Escr}_A$}
			\begin{aligned} 
				\max_{\x \in \mathcal X} & \quad f(\x) \triangleq \mathbb{P}\left\{ \zeta\mid \zeta \in \Kscr_{\Escr}, |\zeta^{\intercal} \x| \leq 1\right\}
			\end{aligned}  
		\end{align}
		if and only if \uss{$\x$ is a solution} of $($\ref{main_prob}$^2_{A})$, defined as
		\begin{align} \tag{\ref{main_prob}$^{2}_{A}$} 
			\uss{\max_{\x\in \mathcal{X}}} & \quad \uss{g(\x)}\triangleq \mathbb{P}\left\{ \ib{\eta}\mid \ib{\eta} \in \Kscr_{2}, |\eta^{\intercal} \uss{\Sigma^{1/2}U \x}| \leq 1\right\}.     
		\end{align}
	\end{lemma}  
	\uss{We may then solve (\ref{main_prob}$^{2}_{A}$) by adding a new variable $\y$, defined using the linear constraint $\y = \Sigma^{1/2} U\x$, and  restate (\ref{main_prob}$^{2}_{A}$) as follows.
		\begin{align} \label{ext-prob}
			\max_{\x\in \mathcal{X}, \y = \Sigma^{1/2} U \x} &  \mathbb{P}\left\{ \ib{\eta}\mid \ib{\eta} \in \Kscr_{2}, |\eta^{\intercal} \uss{\y}| \leq 1\right\}.
		\end{align}
		By our representation theorem  (Theorem~\ref{theo-rep}), 
		we may then claim the following from Lemma~\ref{equiv-ellipsoid} (proof omitted).} 
	\begin{proposition}[Representation of (\ref{main_prob}$^{\Escr}_A$) via Gaussian transformation]
		Consider the function $f$ in \eqref{main_prob}. Suppose Assumption~\ref{dist-ass} holds and $c(\x,\zeta)$ is defined as (\ref{type1}).  Suppose $\Kscr$ is an ellipsoid in $\Real^n$, \afj{defined in \eqref{ellipsoid}}.
		\uss{Suppose $(\x,\y)$ is a solution of \eqref{ext-prob}. Then  $\x$ is  a solution to $(\ref{main_prob}_{A}^\Escr)$.} 
	\end{proposition}

	\subsection{Expectation-valued convex representations for Setting (B)}\label{sec:2.3}
	{{In this section, we consider the regime where
			$c(\x,\zeta) = T \x - \zeta$, $T \in \Real^{d \times n}$, and $\zeta \in \Real^d$. We denote the $i$th component of $c$ by $c_i$ while the $i$th row of $T$ is denoted by $T_{i, \bullet}$. 
			The following proposition articulates a convex counterpart of
			(\ref{main_prob}).  
			
			\begin{theorem}[Transformation of (\ref{main_prob}$_B$) to convex program] Consider the problem \eqref{main_prob}. Suppose Assumptions~\ref{dist-ass} holds and $c(\x,\zeta)$ is
				defined as \eqref{type2}. Suppose $h(\x)
				\triangleq -\log\left( f(\x)\right).$ Then the following hold. 
				
				\begin{enumerate}
					\item[(a)] $h$ is convex on $\mathcal X$. 
					\item[(b)] $\x^*$ minimizes $h$ over $\mathcal{X}$ if and only if $\x^*$ is a maximizer of $f$ over $\mathcal{X}$.
				\end{enumerate}
			\end{theorem}
			\begin{proof}
				
				\noindent (a) We observe that $f(\x) = \mathbb{P}\left(\zeta \mid T\x \geq \zeta \right)
				=H_{\zeta}(T\x)$, where $ H_{\zeta} $ is the probability distribution function of
				the random vector $ \zeta $.  If $ T\x \in \mathbb{R}^d$, $ H_{\zeta}$ is a log-concave function~\cite[Theorem
				4.2.4.]{prekopa2013}, implying that $
				\log(f) $ is a concave function. Hence, it follows that $ h$ is a convex
				function where $h(\x) =-\mbox{log}(f(\x))$.
				
				\noindent (b) Suppose $\x^*$ is a maximizer of $f$ over $\mathcal X$. Then $ f(\x) \leq f(\x^*) $ for all $ \x \in \mathcal X $. Since $ h = -\log(f) $ is a monotonically decreasing function for $ f(\x) > 0 $ we have that $f(\x) \leq f(\x^*)$ if and only if $-\log(f(\x)) \geq -\log(f(\x^*))$ for all $\x \in \mathcal X.$  Consequently, $\x^*$ is a minimizer of $h$ over $\mathcal X$. 
			\end{proof}
			
			{
				\begin{proposition}[Representation of $f$ \ieb{as expectation for symmetric $\Kscr$}] \label{prop:setting B}
					Consider the problem \eqref{main_prob}. Suppose Assumptions~\ref{dist-ass} holds, $c(\x,\zeta)$ is defined as \eqref{type2}, and for all $\x \in \mathcal X$, $T_{i,\bullet}\x \geq \delta$ for every $i$ and for some $\delta > 0$. Then the following equality holds.
					\begin{align}
						\label{type2-exp} & \ \mathbb{P} \left\{\zeta: \zeta \in \K(\x)\right\} = \int_{\mathbb{R}^\ieb{d}} F(\x,\xi) \tilde p(\xi) d\xi = \mathbb{E}_{\tilde p(\xi)} [F(\x,\xi)], \ \mbox{ where } \\ 
						\notag F(\x,\xi) & \triangleq  \mathcal{C}_{\Kscr}(2\pi)^{\ieb{d}/2} 
						e^{-g(\ieb{\x},\xi)+\frac{\|\xi\|^2_{\Kscr}}{2}}, 
						\mathcal{C}  \triangleq \tfrac{1}{\mathrm{Vol}(\Kscr)} \ieb{\tfrac{1}{\Gamma(1+d/2)}},
						\tilde p(\xi) \triangleq 
						\tfrac{1}{ (2\pi)^{\ieb{d}/2} D_{\Kscr} } e^{\frac{-\|\xi\|^2_{\Kscr}}{2}},\\
						&\mbox{ where } \notag g(\ieb{\x}, \xi) \triangleq \max\left\{ \left(\frac{\max(\xi_1,0)}{T_{1,\bullet}\x}\right)^2,\cdots, \left(\frac{\max(\xi_d,0)}{T_{d,\bullet}\x}\right)^2,\|\xi\|_{\Kscr}^2  \right\}, \mbox{ and }
					\end{align} 
					$D_{\Kscr}$ is a positive scalar such that $\int_{\Real^\ieb{d}}  
					\tfrac{1}{ (2\pi)^{\ieb{d}/2} D_{\Kscr} } e^{\frac{-\|\xi\|^2_{\Kscr}}{2}} = 1$, and $\mathcal{C}_{\Kscr} \triangleq \ieb{\mathcal{C}}D_{\Kscr}$.
				\end{proposition}
				\begin{proof}
					In this instance, the set $ \K(\x) $ is defined as following:
					\begin{align*}
						\mathbf{K}(\x)  &  \triangleq \left\{ \zeta: \zeta \in \Kscr\}\cap \{\zeta: \zeta \leq T\x  \right\}.
					\end{align*}
					Since the set $\mathcal{K}$ is convex, compact, and symmetric, the Minkowski functional of
					$\mathcal{K}$ defines a norm, and hence, it is a PHF. Moreover, by
					the definition of the Minkowski functional we have
					$\zeta \in \mathcal{K} \Leftrightarrow  \|\zeta\|_\mathcal{K} \leq 1.$
					By using this definition, we may rewrite $\K(x)$ as follows,
					where by assumption $T_{i,\bullet} \x \geq
					\delta > 0$ for all $\x$ implying that $T_{i,\bullet} \x  >  0$ for $i=
					1, \cdots, d$. 
					\begin{align*}
						\K(\x) &= \{\zeta: \|\zeta\|_\mathcal{K} \leq 1 \} \bigcap \left\{\zeta: \bigcap_{i=1}^d \frac{\max\uss{\{\zeta_i,0\}}}{T_{i,\bullet} \x}\leq 1  \right\}\\
						&= \{\zeta: \|\zeta\|^{2}_\mathcal{K} \leq 1 \} \bigcap \left\{\zeta: \bigcap_{i=1}^d \left(\frac{\max\uss{\{\zeta_i,0\}}}{T_{i,\bullet} \x}\right)^{2} \leq 1  \right\}\\
						&= \uss{\left\{\zeta: \max\left\{ \|\zeta\|_{\Kscr}^2,\left(\frac{\max\{\zeta_1,0\}}{T_{1,\bullet}\x}\right)^2,\cdots, \left(\frac{\max\{\zeta_d,0\}}{T_{d,\bullet}\x}\right)^2  \right\}\leq 1 \right\}},
					\end{align*}
					where the squared expression is employed to obtain a PHF of degree $2$. 
					Since $ \uss{g_i(\x,\zeta) \ \triangleq  \ } \uss{\left(\frac{\max\{\zeta_i,0\}}{T_{i,\bullet} \x}\right)^2} $ for $ i=1,\ldots, d $ and $ g_{d+1}(\uss{\x}, \zeta) \triangleq \|\zeta\|_{\Kscr}^{\uss{2}} $ are PHFs of degree $ 2$, then $ g(\uss{\x}, \zeta) \triangleq \max \{ g_1(\uss{\x}, \zeta),\ldots,g_{d+1}(\uss{\x}, \zeta) \} $ is positively homogeneous of degree $ 2 $. By selecting $h(\zeta) = 1$ and $\ieb{\Lambda} = \K(\x)$, we may invoke Lemma~\ref{las-1}, leading to the following equality. 
					\begin{align} \label{eq:f_2}
						f(\x) = 	\frac{1}{\mbox{Vol}(\Kscr)} \int_{\mathbf{K}(\x)} 1 \ d\xi = 
						\frac{1}{\mbox{Vol}(\Kscr)}\frac{1}{\Gamma(1+n/2)}\int_{\mathbb{R}^{\ieb{d}}} e^{-g(\uss{\x},\xi)} d\xi. 
					\end{align}
					In fact the expression \eqref{eq:f_2} can be restated as 
					\begin{align}
						f(\x)  \ & =  \int_{\mathbb{R}^\ieb{d}} \underbrace{\left(\mathcal{C}_{\Kscr} (2\pi)^{\ieb{d}/2} e^{-g(\uss{\x},\xi) +\frac{\|\xi\|^2_{\Kscr}}{2}}\right)}_{\triangleq F(\x,\xi)}\underbrace{\left( \tfrac{1}{D_{\Kscr}(2\pi)^{\ieb{d}/2}} {e^{-\tfrac{\|\xi\|^2_{\Kscr}}{\ieb{2}}}}
							\right)}_{\triangleq \tilde p(\xi)}  \ d\xi \nonumber \\ 
						\notag & =   \int_{\mathbb{R}^\ieb{d}} F(\x,\xi) \ \tilde p(\xi) \ d\xi  = \mathcal{C} \ \mathbb{E}_{\tilde p(\xi)}[F(\x,\xi)], 
						\mbox{ where } \mathcal{C}  \triangleq \tfrac{1}{\text{Vol}(\Kscr)} \ \ieb{\tfrac{1}{\Gamma(1+d/2)}}, 
					\end{align} 
					$\tilde p(\xi)$  denotes the density, 
					$D_{\Kscr}$ is such that $\int_{\Real^\ieb{d}} \tilde{p}(\xi) d\xi = 1$, and $\mathcal{C}_{\Kscr} \triangleq \mathcal{C} D_{\Kscr}$. 
				\end{proof}
				
				\uss{Next, we provide an analog of Lemma~\ref{lemma_Clarke_F} for this integrand $F(\x,\xi)$ but omit the proof. Note that $g_i(\bullet,\xi)$ is differentiable for every $\xi$ and $i = 1, \cdots, d+1$. }
				\uss{\begin{lemma}\label{lemma_Clarke_F2}
						Consider the function $F(\bullet,\xi)$ defined as 
						
						\begin{align*}
							F(\x,\xi) = \begin{cases} \left(\mathcal{C}_{\Kscr} (2\pi)^{\ieb{d}/2} e^{-g_i(\uss{\x},\xi) +\frac{\|\xi\|^2_{\Kscr}}{2}}\right),   & \xi \in \Xi_i(\x)  \\ 
								\left(\mathcal{C}_{\Kscr} (2\pi)^{\ieb{d}/2} e^{-g(\uss{\x},\xi) +\frac{\|\xi\|^2_{\Kscr}}{2}}\right),   & \xi \in \Xi_0(\x), \ \mbox{ where }
							\end{cases}
						\end{align*}
						
						\begin{align*}
							\Xi_i(\x) &\triangleq \left\{\xi \mid  g_i(\x,\xi) > g_j(\x,\xi) \quad  \forall j \neq i \right\}, i = 1, \cdots, d+1 \\ 
							\Xi_0(\x) &\triangleq \left\{\xi \mid  g_i(\x,\xi) = g_j(\x,\xi), \quad \forall (i,j) \in {\cal I}, \right. \\
							& \left. \qquad \quad  \ieb{g_i(\x,\xi)>g_l(\x,\xi) \quad \forall l \notin {\cal I},} \quad {\cal I} \subseteq \{1, \cdots, d+1\}  \right\}.
						\end{align*}
						Then the following hold.
						
						\noindent (a)  $F(\bullet,\xi)$ is locally Lipschitz  for every $\xi$. 
						
						\noindent (b) $F(\bullet,\xi)$ is a Clarke regular function for almost every $\xi \in \Real^\ieb{d}$. 
						
						\noindent (c) For any $\x \in \Real^n$, $ \partial \mathbb{E}[F(\x,\xi)] = \mathbb{E}[\partial F(\x,\xi)].$    
				\end{lemma}}
				
				We now analyze $F(\x,\xi)$ and $G(\x,\xi)$ where $G(\x,\xi) \in \partial F(\x,\xi)$.
				\begin{lemma}[\ieb{Properties of $F(\x,\xi)$ and $G(\x,\xi)$ under symmetric $\Kscr$}] \label{bound sub settingB}
					Consider the problem \eqref{main_prob}. Suppose Assumption~\ref{dist-ass} holds, $c(\x,\zeta)$ is defined as \eqref{type2}, and  {for all $\x \in \mathcal X$, $T_{i,\bullet}\x \geq \delta$ for every $i$} for some $\delta > 0$. Suppose $f(\x)$ is defined as in Prop.~\ref{prop:setting B}. Then the following hold. 
					
					\noindent (a) For any $\x \in \mathcal{X}$, $|F(\x,\xi)|^2 \leq \mathcal{C}^2_{\Kscr} (2\pi)^\ieb{d}$ for every $\xi \in \Real^d$. 
					
					\noindent (b) Given an $\x \in \mathcal{X}$ such that $T_{i,\bullet}\x \geq \delta$ for every $i$ \afj{and} for some $\delta > 0$, and $G(\x,\xi) \in \partial F(\x,\xi)$,  \afj{then it holds that}
					$
					\mathbb{E}[\|G(\x,\xi)\|^2] \leq  16\mathcal{C}_{\Kscr}^2(2\pi)^\ieb{d}  \frac{\|T_{i,\bullet}\|^2}{\delta^2 e^2}
					$
				\end{lemma}
				
				\uss{\begin{proof}
						\uss{Recall the definition of $F(\x,\xi)$ from the statement of Lemma~\ref{lemma_Clarke_F2}. Suppose $\Xi_{i,0}(\x) \triangleq \{ \xi: g_i(\x,\xi) \geq g_j(\x,\xi), j \neq i\}$ for $i = 1, \hdots, d+1$. It can be seen that $\cup_{i=1}^{d+1} \Xi_{i,0}(\x) = \Real^\ieb{d}$.
							We prove (a) by considering two cases. Case (i): $\xi \in \Xi_{i,0}(\x)$ for $i = 1, \cdots, d$.  It follows that 
							\begin{align*} |F(\x,\xi)|^2 = \mathcal{C}^2_{\Kscr} 
								\left( (2\pi)^{\ieb{d}} e^{-2g(\x,\xi)+\|\xi\|_{\Kscr}^2}  \right) 
								\leq \mathcal{C}^2_{\Kscr} \left( (2\pi)^{\ieb{d}} e^{-2g_i(\x,\xi)+ g_i(\x,\xi)}  \right)  \leq \mathcal{C}^2_{\Kscr} (2\pi)^\ieb{d}.
							\end{align*}
							Case (ii): $\xi \in \Xi_{d+1,0}(\x).$ Proceeding similarly, we obtain that
							\begin{align*}
								|F(\x,\xi)|^2
								& \leq \mathcal{C}^2_{\Kscr} \left( (2 \pi)^{\ieb{d}} e^{\us{-2g(\x,\xi)}+{\us{\|\xi\|^2_{\Kscr}}} } \right) \leq  
								\mathcal{C}^2_{\Kscr} \left( (2 \pi)^{\ieb{d}} e^{\us{-2\|\xi\|_{\Kscr}^2}+{\us{\|\xi\|^2_{\Kscr}}} } \right) \leq \mathcal{C}^2_{\Kscr} (2\pi)^\ieb{d}.
							\end{align*}
							Consequently, $|F(\x,\xi)|^2 \leq \mathcal{C}^2_{\Kscr} (2\pi)^\ieb{d}$ for $\xi \in \Real^\ieb{d}$.}  
						
						\noindent (b). We observe that $\partial F(\x,\xi)$ is defined as follows.
						\begin{align*} \scriptsize
							\partial F(\x,\xi) = \begin{cases}  \left( \ieb{\mathcal{C}_{\Kscr}}(2\pi)^{\ieb{d}/2} \frac{2(\max\{\xi_i,0\})^2 T_{i,\bullet}^T}{(T_{i,\bullet} \x)^3}  e^{-g_i(\x,\xi)+\frac{\us{\|\xi\|_{\Kscr}^2}}{2}}\right),  & \xi \in \Xi_i(\x), i = 1, \hdots, d \\ 
								\left( -\ieb{\mathcal{C}_{\Kscr}}(2\pi)^{\ieb{d}/2}  e^{-g(\x,\xi)+\frac{\us{\|\xi\|_{\Kscr}^2}}{2}}\right) H(\x,\xi),  & \xi \in \Xi_0(\x)\\
								{\bf 0}.  & \xi \in \Xi_{d+1}(\x), 
							\end{cases}
						\end{align*}  
						where $H(\x,\xi)$ denotes the Clarke generalized gradient of $g(\x,\xi)$, defined as  
						\begin{align} \label{eq-clarke-H}
							H(\x,\xi) = \left\{ \sum_{\ell \in  {\cal I}} \alpha_\ell \beta_\ell \mid \alpha_\ell \geq 0, \sum_{\ell\in {\cal I}} \alpha_\ell = 1, \beta_{\ell} = \nabla_{\x} g_{\ell}(\x,\xi) \right\}. 
						\end{align}
						\noindent {Proceeding as in Prop~\ref{bound sub}, we} have that 
						$\mathbb{E}_{\tilde{p}} \left[\|G(\x,\xi)\|^2\right]$ can be expressed as follows.
						\begin{align}\notag
							& \quad \mathbb{E}_{\tilde{p}}\left[\|G(\x,\xi)\|^2\right]  = \int_{\Real^n} \|G(\x,\xi)\|^2 \tilde{p}(\xi) d\xi \\
							& = \sum_{i=1}^{d} \int_{\Xi_i(\x)} \|G(\x,\xi)\|^2 \tilde{p}(\xi) d\xi \notag 	
							+ \int_{\Xi_{d+1}(\x)} \|\underbrace{G(\x,\xi)}_{ \ = \ {\bf 0}}\|^2 \tilde{p}(\xi) d\xi 	
							\\ & + \int_{\Xi_0(\x)} \|G(\x,\xi)\|^2 \tilde{p}(\xi) d\xi \label{decomp-int3}\\
							\notag			&  = \sum_{i=1}^{d}\int_{\Xi_i(\x)} \|G(\x,\xi)\|^2 \tilde{p}(\xi) d\xi, 	
						\end{align} 
						where the last equality follows from observing that $G(\x,\xi) = 0$ for $\xi \in \Xi_{d+1}(\x)$ and the integral in \eqref{decomp-int3} is zero because $\Xi_0(\x)$ is a measure zero set. 
						It follows that $\mathbb{E}[\|G(\x,\xi)\|^2]$ can be bounded as follows:
						\begin{align*}
							\mathbb{E}&[\|G(x,\xi)\|^2] \\
							&= \sum_{i=1}^{d} \int_{\Xi_i(x)} 4\mathcal{C}_{\Kscr}^2(2\pi)^\ieb{d} \left(\frac{\xi_{i}^2 T^{\intercal}_{i,\bullet}}{(T_{i,\bullet} \x)^3}\right)^T\left(\frac{\xi_{i}^2 T^{\intercal}_{i,\bullet}}{(T_{i\bullet} \x)^3}\right) e^{-\frac{2(\xi_i)^2}{(T_{i,\bullet} \x)^2}+\|\xi\|_{\Kscr}^2}\left( \tfrac{1}{D_{\Kscr}(2\pi)^{\ieb{d}/2}} e^{-\tfrac{\|\xi\|_{\Kscr}^2}{2}}\right)d\xi,\\
							&=\sum_{i=1}^{d} \int_{\Xi_i(x)} 4\mathcal{C}_{\Kscr}^2(2\pi)^\ieb{d} \frac{\|T_{i,\bullet}\|^2}{(T_{i,\bullet} \x)^2}\left(\frac{\xi_i}{{T_{i,\bullet} \x}}\right)^4 e^{-\frac{2(\xi_i)^2}{(T_{i,\bullet} \x)^2}+\|\xi\|_{\Kscr}^2}\left( \tfrac{1}{D_{\Kscr}(2\pi)^{\ieb{d}/2}} e^{-\tfrac{\|\xi\|_{\Kscr}^2}{2}}\right)d\xi,\\
							&\leq \sum_{i=1}^d \int_{\Xi_i(x)} 4\mathcal{C}_{\Kscr}^2(2\pi)^\ieb{d} \frac{\|T_{i,\bullet}\|^2}{\delta^2}\left(\frac{\xi_i}{{T_{i,\bullet} \x}}\right)^4 e^{-\frac{(\xi_i)^2}{(T_{i,\bullet} \x)^2}}\left( \tfrac{1}{D_{\Kscr}(2\pi)^{\ieb{d}/2}} e^{-\tfrac{\|\xi\|_{\Kscr}^2}{2}}\right)d\xi\\
						\end{align*}
						where the inequality follows from $\ieb{\xi \in \Xi_i(\x) \mbox{ for } i = 1, \hdots, d }$ and $ T_{i,\bullet} \x \geq \delta $ for all $ i $. Next, we consider the expression $\left(\frac{\xi_i}{{T_{i,\bullet} \x}}\right)^4 e^{-\left(\frac{\xi_i}{{T_{i,\bullet} \x}}\right)^2}$ or $u^2 e^{-u}$. We note that by Lemma~\ref{unimodal}, $u^* = 2$ is a maximizer with value $4e^{-2}$. Hence, we have that 
						\begin{align*}
							\mathbb{E}[\|G(\x,\xi)\|^2] &\leq  \sum_{i=1}^d \int_{\Xi_i(x)} 4\mathcal{C}_{\Kscr}^2(2\pi)^\ieb{d} \frac{\|T_{i,\bullet}\|^2}{\delta^2}\frac{4}{e^2}\left( \frac{1}{D_{\Kscr}(2\pi)^{\ieb{d}/2}} e^{-\tfrac{\|\xi\|_{\Kscr}^2}{2}}\right)d\xi\\
							&= 16\mathcal{C}_{\Kscr}^2(2\pi)^\ieb{d}   \sum_{i=1}^d \frac{\|T_{i,\bullet}\|^2}{\delta^2 e^2} \int_{\Xi_{i}(\x)}\tfrac{1}{D_{\Kscr}(2\pi)^{\ieb{d}/2}} e^{-\tfrac{\|\xi\|_{\Kscr}^2}{2}} d\xi \\
							&   =  
							16\mathcal{C}_{\Kscr}^2(2\pi)^\ieb{d} \sum_{i=1}^d \frac{\|T_{i,\bullet}\|^2}{\delta^2 e^2}. 
						\end{align*} 
				\end{proof}}
				
				\ieb{We now specialize these results to regimes where $\Kscr$ is an $\ell_p$-ball in $\Real^n$ and not necessarily symmetric about the origin via the following Proposition.}
				
				\begin{assumption}\label{dist-ass-3}
					\ieb{The random variable $\zeta$ is uniformly distributed on the set $\Kscr \subset \Real^d$ where $\Kscr \triangleq \{ \zeta : \ \| \zeta - \mu \|_{p} \leq \alpha\}$. The set $\mathcal X$ is closed, convex, and bounded. }
				\end{assumption}
				
				\begin{proposition}[Representation and boundedness  under asymmetric $\Kscr$]\label{bound-asym-sett_B} 
					Consider the problem \eqref{main_prob}. Suppose \ieb{Assumption~\ref{dist-ass-3}} holds and there exists $\delta > 0$ such that $T_{i,\bullet}\x -\mu_i \geq \delta$ for $i = 1, \cdots, d.$  Then the following hold. 
					
					\noindent (a) 
					$f(\x) \triangleq \mathbb{E}_{\tilde p(\xi)}[F(\x,\xi)]$, where $\sigma^2 = \alpha^2$,\\ 
					$F(\x,\xi)  \triangleq  \mathcal{C}(2\pi\ieb{\sigma^2})^{\ieb{d}/2} 
					e^{-g(\ieb{\x},\xi)+\frac{\|\xi\|^2_{2}}{2\sigma^2}}$, $ \mathcal{C}  \triangleq \tfrac{1}{\mathrm{Vol}(\Kscr)} \ieb{\tfrac{1}{\Gamma(1+d/2)}}$, $ \tilde p(\xi) \triangleq 
					\tfrac{1}{ (2\pi\ieb{\sigma^2})^{\ieb{d}/2}  } e^{\frac{-\|\xi\|^2_{2}}{2\sigma^2}},$
					and,
					$g(\ieb{\x},\xi) \triangleq \max\left\{ \frac{1}{\alpha^2}\|\xi\|_{p}^2,\left(\frac{\max(\xi_1,0)}{T_{1,\bullet}\x-\mu_1}\right)^2,\cdots, \left(\frac{\max(\xi_d,0)}{T_{d,\bullet}\x-\mu_d}\right)^2  \right\}.$
					
					\noindent (b) Given an $\x \in \mathcal X$ and $T\x - \mu \geq \delta e$, $F(\x,\xi) \geq 0$ and $|F(\x,\xi)|^2 \leq  \mathcal{C}^2 (2\pi\sigma^2)^\ieb{d} $ for every $\xi \in  \Real^d.$
					
					\noindent (c) Given an $\x \in \mathcal X$ and $T\x - \mu \geq \delta e$, and $G(\x,\xi) \in \partial F(\x,\xi)$, \ieb{then it holds that}  
					\begin{align}
						\mathbb{E}[\|G(\x,\xi)\|^2] \leq 16\mathcal{C}^2(2\pi\sigma^2)^\ieb{d} \frac{\|T_{i,\bullet}\|^2}{\delta^2 e^2}.
					\end{align}
				\end{proposition}
			}
			\uss{Before concluding, we comment on the assumptions employed in this section. 
				\noindent (a) {\em Assumptions on $\mathbb{\zeta}$.} We assume that $\zeta$
				is  uniformly distributed on the set $\mathcal{K}$ which is a compact and
				convex set, symmetric about the origin. The second requirement is that $\mathcal{K}$ be convex,
				compact, and symmetric (with corollories provided for specializing these
				results to an ellipsoid). This property
				allows us to claim that the Minkowski function of $\mathcal{K}$ is a norm, a
				key step in the analysis. However, we do develop an extension to non-symmetric
				regimes where $\mathcal{K} \triangleq \{\zeta \mid \|\zeta - \mu\|_p \leq
				\alpha\}$. 
				The requirement that $\zeta$ is
				uniformly distributed may be weakened to log-concave measures and this will be
				the focus of future work, as noted in the concluding section.\\
				
				\noindent (b) {\em Definition of $c(\x,\zeta)$ in Settings A and B.} We have adopted two distinct choices for $c$; i.e. $c(\x,\zeta) = 1 - | \zeta^T\x|^m$ (Setting A) and $T\x - \zeta$ (Setting B). Extensions to this are also possible where $\zeta^T A\x$ is employed in Setting A. This can be easily addressed by adding variables. More general extensions will be considered as part of future work.}
			
			\section{An efficient stochastic approximation framework}\label{sec:3} 
			In the prior section, we observed that the function $f$ could be
			recast as an expectation of $F(\x,\xi)$ with respect to a suitable
			density function. In Section~\ref{sec:3.1}, we cast the stochastic
			optimization of problem as a {\em convex compositional stochastic optimization problem} and comment on why available schemes do not suffice. We then provide some background and define
			the algorithmic framework in Section~\ref{sec:3.2}. Finally,
			convergence and rate analysis are provided in Section~\ref{sec:3.3}.
			
			\subsection{Convex compositional stochastic optimization problem} \label{sec:3.1}
			The optimization problem of interest, denoted by {\bf (PM)}, can be cast as the following convex compositional optimization problem. 
			{\begin{align}\label{def-h}
					\min_{\x \in \mathcal{X}} \ h(\x), \mbox{ where } h(\x) \triangleq \psi(\mathbb{E}[F(\x,\xi)])  \mbox{ and } \psi(\y) \triangleq \begin{cases}  \frac{1}{y}  &  (\ref{main_prob}_A)  \\
						\ieb{- \log\left(y\right)}. & (\ref{main_prob}_B) 
					\end{cases}
			\end{align}}
			Before proceeding, we provide a brief review of SA schemes and their variance-reduced and compositional counterparts.
			
			\underline{ (a) \em Stochastic approximation (SA) schemes.} SA schemes
			represent a class of techniques rooted in the seminal work
			by Robbins and Monro~\cite{robbins1951}. In the last several decades, there has been a
			tremendous amount of research in stochastic approximation applied to minimizing
			a convex function $h$, defined as $h(\x) \triangleq \mathbb{E}[H(\x,\omega)]$
			over a closed and convex set $\mathcal X$. Noteworthy amongst these being the
			long-step averaging framework by ~\cite{polyak1990} and ~\cite{polyak1992}.
			In fact, in ~\cite{nemirovski2009} the authors developed a {\em robust} stochastic
			approximation framework for convex stochastic optimization in which a constant
			steplength of prescribed size was employed over a pre-selected number of SA
			steps. Such a scheme admits the optimal rate of convergence of
			$\mathbb{E}[h(\bar{\x}_K) - h^*] \leq \mathcal{O}(1/\sqrt{K})$ where $K$
			represented the number of steps and $\bar{\x}_K$ denotes the iterate average
			over $K$ steps. \\

			\underline{(b) \em Variance-reduced schemes.} A key shortcoming of SA schemes is
			the gap in the convergence rate between the deterministic schemes and their SA
			analogs. \uss{This gap is particularly irksome in the presence of complicated constraints, since the projection operation is computationally expensive and in such cases, deterministic rates of convergence have profound benefits. For instance, to compute an $\epsilon$-solution for a smooth convex expectation-valued problem, traditional SA schemes require at most $\mathcal{O}(1/\epsilon^2)$ while the variance-reduced counterparts require $\mathcal{O}(1/\epsilon)$. For instance, if $\epsilon=1$e-$3$, standard SA schemes require $\mathcal{O}(1$e$6)$ projection steps while variance-reduced counterparts  require $\mathcal{O}(1$e$3)$ steps, a significant difference. When considering sample-complexity, we note that in some instances such as ~\cite{jalilzadeh2018optimal}, one may be able to (nearly) match the sample complexity of $\mathcal{O}(1/\epsilon^2)$. These improved rates are achieved by either utilizing an increasing
				batch-size of gradients or by solving a sequence of stochastic subproblems to
				increasing degrees of inexactness. Such avenues have derived deterministic rates of convergence in 
				smooth strongly convex~\cite{byrd12,shanbhag15budget,xie20siadmm}, smooth convex~\cite{ghadimi2016accelerated}, nonsmooth convex~\cite{jalilzadeh2018optimal}, and nonconvex regimes~\cite{ghadimi2016accelerated,lei20asynchronous}}. 
			Notably, in many of
			these settings, the schemes admit optimal or near-optimal sample
			complexities~\cite{shanbhag15budget,xie20siadmm,jalilzadeh2018optimal,lei20asynchronous}.\\

			\underline{(c) Compositional stochastic optimization.}  \uss{The earliest
				efforts on compositional optimization appear to be the almost-sure convergence
				guarantees provided by Ermoliev~\cite{Erm76} for two-level problems. Rate
				statements~\cite{wang2017stochastic,wang2017accelerating} and
				variance-reduction (in finite sample-space regimes)~\cite{lian2017finite} has
				been studied while  multi-level settings were first considered
				by~\cite{yang2019multilevel}. Optimal sample-complexity in nonconvex regimes
				was shown for two-level~\cite{ghadimi2020single} and
				multi-level~\cite{balasubramanian2020stochastic,chen2021solving} regimes.
				However, when the inner function is nonsmooth (as in this setting), the best
				known rate has been provided in ~\cite{wang2017stochastic} where a rate of
				$\mathcal{O}(k^{-1/4})$ has been derived.} \uss{We note that in the present
				setting, sample complexity is of less relevance since $\xi$ is a Gausian random
				variable and sampling is cheap with no explicit limitations on data (unlike in
				finite-sum machine learning problems). Instead, in this setting, we argue that
				iteration complexity is of more relevance.}\\  
			
			
			\uss{\underline{ (d) \em Gaps and shortcomings in existing SA and compositional SA schemes.}
				A  prototypical SA scheme for minimizing a convex function $h$,
				defined as $h(\x) \triangleq \mathbb{E}[H(\x,\xi)]$ and given $\x_0 \in \mathcal{X}$, generates a sequence $\{\x_k\}$ as
				follows.\begin{align} \x_{k+1} := \Pi_{\mathcal X} \left[ \x_k - \gamma_k
					d(\x_k,\xi_k) \right], \qquad k \geq 0 \end{align} where $d(\x_k,\xi_k)$ is
				assumed to be a sampled subgradient, the interchange between the expectation
				and subdifferential operator is assumed to hold for any $\x$, i.e. $\partial
				\mathbb{E}[H(\x,\xi)] = \mathbb{E}[\partial H(\x,\xi)]$,  and
				$\mathbb{E}[d(\x_k,\xi) \mid \x_k]  \in \partial_{\x}
				\mathbb{E}[H(\x_k,\xi)]$. When contending with $\psi(\mathbb{E}[F(\x,\xi)])$, by the chain rule~\cite{clarke98}, we have that 
				\begin{align}
					\partial_{\x} \psi(\mathbb{E}[F(\x,\xi)]) & = \partial_{\x} [\mathbb{E}[F(\x,\xi)]] \psi'(\mathbb{E}[F(\x,\xi)]) \notag \\
					& =\mathbb{E}[  
					\partial_{\x} F(\x,\xi)]\psi'(\mathbb{E}[F(\x,\xi)]),
				\end{align}
				where  the second equality
				is a consequence of invoking Lemma~\ref{lemma_Clarke_F}. Consequently, an unbiased
				subgradient of $\psi(\mathbb{E}[F(\x,\xi)])$ is given by
				$G(\x,\xi)\psi'(\mathbb{E}[F(\x,\xi)])$ and requires access to
				$\psi'(\mathbb{E}[F(\x,\xi)])$; however, the latter cannot be accessed and
				therefore an unbiased subgradient cannot be tractably evaluated and standard SA
				schemes cannot be adopted.}\\

			\uss{ \underline{(e) Related numerical schemes.} 
				\begin{enumerate}
					\item[(i)] {\em SA and Mini-batch SA schemes.} In the context of stochastic optimization with conditionally unbiased gradients being available,  single-sample SA schemes are characterized by an optimal rate of $\mathcal{O}(k^{-1/2})$ while mini-batch variants employ a gradient estimator with reduced bias. However, in the current regime, such estimators are complicated by bias.  Note that these schemes are not equipped by either asymptotic convergence or non-asymptotic rate guarantees.
					Yet, given that such schemes enjoy an optimal rate in unbiased regimes, SA schemes and mini-batch variants of SA provide a useful benchmark of comparison.
					\item[(ii)] {\em Compositional SA schemes.}  The presence of bias arising from the presence of a compositional structure has been addressed by compositional stochastic approximation schemes by adding a parallel updating scheme~\cite{wang2017stochastic}. In nonsmooth regimes, such an avenue is characterized by a convergence rate of $\mathcal{O}(k^{-1/4})$ while our proposed scheme achieves a rate of approximately $\mathcal{O}(k^{-1/2})$. Note that sample (or oracle) complexity is less relevant here since sampling is (relatively) cheap and data is not limited by any means. Given the significant difference in rates, we have not provided an additional comparison with compositional SA schemes in the current manuscript.
					\item[(iii)] {\em Sample-average approximation via integer programming.} Finally, the other competing approach for computing global minimizers of chance-constrained problems is via sample-average approximation (SAA) where the SAA problem is resolved via integer programming~\cite{ahmed09sample}. We introduce this comparison to demonstrate the difference in scalability and from the standpoint that  this avenue also provides an additional certification that our proposed (r-VRSA) scheme is indeed finding near-global minimizers.
			\end{enumerate}}

			\subsection{Background and Algorithm definition}\label{sec:3.2} 
			
			\uss{We observe that the problem of interest is $\min_{\x \in \Xscr}
				\psi(\mathbb{E}[F(\x,\xi)])$. We first provide a result that allows us to claim
				that $\partial_{\x} [\psi(\mathbb{E}[F(\x,\xi)])] = 
				\psi'(\mathbb{E}[F(\x,\xi)] \partial_{\x} [F(\x,\xi)]$ indeed holds.
				
				\begin{lemma}
					Suppose $F(\bullet,\xi)$ is a Clarke regular function for every $\xi \in \Xi$, $\psi$ is a continuously differentiable function, and $\Xscr$ is a nonempty, compact, and convex set in $\Real^n$.  Then the following hold.
					
					\noindent (a) $F(\bullet,\xi)$ is Lipschitz continuous on $\Xscr$ with a Lipschitz constant $L(\xi)$ where $\mathbb{E}[L(\xi)] \leq \tilde{L}$.  
					
					\noindent (c) Suppose $f(\x) = \mathbb{E}[F(\x,\xi)]$ and $f(\bar{\x})$ is finite for some $\bar \x \in \Xscr$. Then $f$ is Clarke-regular and for any $\x \in \Xscr$, $\partial_{\x} f(\x) = \mathbb{E}[\partial_{\x} F(\x,\xi)]$. 
					
					\noindent (d) Suppose $\psi: \Real^+ \to \Real$ is a continuously differentiable function. Then for any  $\x \in \Xscr$ such that $f(\x) = \mathbb{E}[F(\x,\xi)] > 0$,  $\partial_{\x} [\psi(f(\x)] = \psi'(\mathbb{E}[F(\x,\xi)]) \mathbb{E}[\partial_{\x} [F(\x,\xi)]]$.  
				\end{lemma}
				\begin{proof}
					\noindent (a)  We present the proof for Setting A. We observe that $F(\bullet,\xi)$ is a piecewise-smooth function for every $\xi$. Then for any $\x \in \Xscr$,  $\|\nabla_{\x} F(\x,\xi)\|$ is bounded as follows at points where $\x$ is smooth:
					\begin{align*}
						\|\nabla_{\x} F(\x,\xi)\| & = \| C_{\Kscr} (2\pi)^{n/2} (2\xi \x^T\x) e^{-\max\{\xi^Tx,\|\xi\|_{\Kscr}^2\} + \|\xi\|^2_{\Kscr}/2} \| \\
						& \leq  C_{\Kscr} (2\pi)^{n/2} (\|\xi\|^2 + (\xi^T\x)^2).
					\end{align*}  
					Consequently, $F(\bullet,\xi)$ is a Lipschitz continuous function with $L(\xi) = C_{\Kscr} (2\pi)^{n/2} (\|\xi\|^2 + (\xi^T\x)^2)$. Furthermore, we have that $\mathbb{E}[L(\xi)] \leq C_{\Kscr} (2\pi)^{n/2} (\mathbb{E}[\|\xi\|^2] + \frac{(\mathbb{E}[\|\xi\|^2] + \|\x\|^2])}{2} < \tilde{L},$ a consequence of the boundedness of the second moment and the compactness of $\mathcal X$.\\  
					
					\noindent (b)  By definition, $f(\x) = \mathbb{E}[F(\x,\xi)]$. Therefore, $f$ is Lipschitz continuous with constant $\tilde{L}$ by utilizing convexity of the norm and Jensen's inequality as well as part (a). \\ 
					
					\noindent (c) Since $F(\bullet,\xi)$ is Clarke-regular on $\Xscr$, $f$ is Lipschitz continuous on $\Xscr$, $f$ is defined at some $\bar \x \in \Xscr$, we have that $f$ is Clarke regular on $\Xscr$ and for any $\x \in \Xscr$, $\partial_{\x} f(\x) = \mathbb{E}[\partial_{\x} F(\x,\xi)]$~\cite[Th.~2.7.2]{clarke98}.\\ 
					
					\noindent (d) This follows from noting that by recalling that $\psi: \Real_+ \to \Real$ is continuously differentiable and $f$ is Lipschitz continuous on $\Xscr$, and then invoking~\cite[Cor.~2.6.6]{clarke98}, we have that     
					$\partial_{\x} [\psi(\mathbb{E}[F(\x,\xi)]] = \psi'(\mathbb{E}[F(\x,\xi)]) \partial_{\x} [\mathbb{E}[F(\x,\xi)]]$.
					
				\end{proof}
				
				Consequently, an unbiased stochastic subgradient of $h$ is given by a
				measurable selection \\ {$h'(\mathbb{E}[F(\x,\xi)] G(\x,\xi)$} where $G(\x,\xi) \in
				\partial_{\x} F(\x,\xi).$ However, such a selection cannot be efficiently
				evaluated since it requires $\mathbb{E}[F(\x,\xi)]$ which is unavailable.
				Instead, we employ a biased variance-reduced counterpart given by the
				following.  \begin{align}  
					D_k = \psi'_{\epsilon}\left(\tfrac{\sum_{j=1}^{N_k} F(\x_k,\xi_{k,j})}{N_k}\right) \tfrac{\sum_{j=1}^{N_k} G(\x_k,\xi_{j,k})}{N_k}.
				\end{align} 
				We observe that the bias in defining the estimator $D_k$ arises from approximating
				$\psi'(\mathbb{E}[F(\x_k,\xi)])$ by $\psi'_{\epsilon}\left(\tfrac{\sum_{j=1}^{N_k}
					F(\x_k,\xi_{k,j})}{N_k}\right)$ where $\psi'_{\epsilon}$ is suitably defined
				approximation of $\psi$ with parameter $\epsilon$. 
				Consequently, we propose the following regularized variance-reduced
				stochastic subgradient scheme for minimizing $h(\x)$ in either
				(\ref{main_prob}$_A$) or (\ref{main_prob}$_B$)}. We define the variance-reduced
			sampled gradient \afj{$D_k$} as follows for each of these settings.
			\begin{align}\label{def-Dk}
				D_k \triangleq  \begin{cases}
					\frac{-(G_k+\bar{w}_{G,k})}{(f(\x_k) + \bar{w}_{f,k})^2+\epsilon_k}, &  (\mbox{Setting A})\\
					{-\frac{(G_{j}+\bar{w}_{G,k})}{(f(\x_k) + \bar{w}_{f,j,k})+\epsilon_k}},&  (\mbox{Setting B}) 
				\end{cases}
			\end{align}
			where $\bar{w}_k \triangleq d_k - D_k$, $d_k \in \partial_{\x} h(\x_k)$, 
			$\bar{w}_{f,k}$ and $\bar{w}_{G,k}$ are defined as   
			\begin{align}\label{def-w}
				\bar{w}_{f,k} & \triangleq \tfrac{\sum_{j=1}^{N_k} F(\x_k,\xi) - f(\x_k)}{N_k} \mbox{ and }   
				\bar{w}_{G,k}  \triangleq \tfrac{\sum_{j=1}^{N_k} G(\x_k,\xi_j) - \mathbb{E}[G(\x_k,\xi)]}{N_k},
			\end{align} respectively.
			We begin by assuming the existence of the following stochastic oracles, crucial for the development of the proposed first-order schemes.  
			\begin{assumption}[Stochastic zeroth and first-order oracles]\label{sfo}
				There exist a stochastic zeroth-order oracle and a stochastic first-order oracle that given $\x$,  
				\uss{produce independent samples $F(\x,\xi)$ and $G(\x,\xi) \in \partial F(\x,\xi)$ in Settings A and B.} 
			
		\end{assumption} 
		We now define the $\sigma$-algebra $\mathcal F_k$ for Setting B (Setting $A$ is defined analogously).
		\begin{align}
			\mathcal{F}_{f,k} & \triangleq \left\{ \{F(\x_0,\xi_j)_{j=1}^{N_0}, \{F(\x_1,\xi_j)_{j=1}^{N_1}, \cdots, \{F(\x_k,\xi_j)_{j=1}^{N_k}\right\},  \\
			\mathcal{F}_{G,k} & \triangleq \left\{ \{G(\x_0,\xi_j)_{j=1}^{N_0}, \{G(\x_1,\xi_j)_{j=1}^{N_1}, \cdots, \{G(\x_k,\xi_j)_{j=1}^{N_k}\right\},  \\
			\mathcal{F}_k & \triangleq \mathcal{F}_{f,k} \cup \mathcal{F}_{G,k} \cup \{\x_0\}.
		\end{align}  
		
		{Suppose $\{\x_k\}$ is a sequence in $\mathcal X$. Then the following result holds.
			\begin{lemma}\label{cond-sec-mom}
				For any $\x_k \in \mathcal X$,  suppose $\bar{w}_{f,k}$ and $\bar{w}_{G,k}$ are defined as in \eqref{def-w}.    
				Then for all $k \geq 0$, \uss{ $\mathbb{E}[\|\bar{w}_{f,k}\|^2 \mid \mathcal{F}_k ] \leq \tfrac{\nu^2_f}{N_k}$ and $\mathbb{E}[\|\bar{w}_{G,k}\|^2 \mid \mathcal{F}_k ] \leq \tfrac{\nu^2_G}{N_k}$ where }
				
				\begin{align} \quad
					\tag{Setting A} & \quad \nu^2_f  \triangleq {2(\mathcal{C}_{\Kscr}^2(2\pi)^n+1)}, 
					\nu^2_G  \triangleq \frac{\mathcal{C}^2_{\Kscr}(2\pi)^n\mathbb{E}_{\tilde p}[\|\xi\|^2]}{e},\\  
					\tag{Setting B} & \quad
					\nu^2_{f}  \triangleq \uss{\mathcal{C}^2 (2\pi\sigma^2)^d} \mbox{ and }
					\uss{\nu^2_{G} \triangleq	 
						16\mathcal{C}_{\Kscr}^2(2\pi)^\ieb{d} \sum_{i=1}^d \frac{\|T_{i,\bullet}\|^2}{\delta^2 e^2}}. 
				\end{align}
			\end{lemma}
		}
		
		\begin{algorithm}[htbp] \caption{\bf {Regularized VR}  stochastic approximation {\bf (r-VRSA)}}
			\em
			\label{smooth scheme}
			(0) Given $\x_0 \in \mathcal X$ and  positive
			sequences $\{\gamma_k,{\epsilon_k}, N_k\}$;  set $k := 1$. \\
			(1) $\x_ {k+1}:=\Pi_{\mathcal X} \left[\x_k- \gamma_k  D_k\right],  \mbox{ where } D_k$ is defined in \eqref{def-Dk} \\
			(2) If $k>K$, then stop; else $k := k+1$; return
			to (1).
		\end{algorithm}
	\begin{assumption}\label{assump_sub}
		There exists an $\epsilon_f$ such that $f(\x_k) \geq \epsilon_f$ and for any $\x_k \in \mathcal X$. For any $\x, \y \in \mathcal X$, $ \|\x-\y\|^2 \leq B^2.$ 
	\end{assumption}
	
	\begin{lemma}\label{bd-sec-moment} Suppose Assumptions~\ref{sfo} and~\ref{assump_sub} hold.  Consider any $\x_k \in \mathcal X$. Suppose $N_k \in \mathbb{Z}_+$ and $\epsilon_k \triangleq \tfrac{1}{N_k^{1/4}}$. Suppose $\bar{w}_k \triangleq D_k - d_k$, where $D_k$ is defined in \eqref{def-Dk} and $d_k \in \partial h(\x_k)$. \ieb{Suppose $ \mathbb{E}[\| G(\x_k,\xi)\|^2 \mid \mathcal{F}_k] \leq M_G^2 $ and $ | F(\x,\xi)| \leq \ib{M_F}$ for any $\x, \xi.$} Then $\mathbb{E}[\|\bar{w}_k\|^2] \leq \tfrac{\nu^2}{\sqrt{N_k}}$, where $\nu^2 = \sum_{j=1}^d \nu_j^2$ in Setting B,  
		\begin{align} 
			\nu^2 & \triangleq \frac{3\nu_G^2}{\epsilon_f^2}+M_G^2  \frac{{24}\nu_f^2}{\epsilon_f^{{4}} }+\frac{{6}({M_F}^2+1)\nu_f^2}{\epsilon_f^{{4}}}  \tag{ Setting A}\\
			\nu^2 & \triangleq 
			\uss{\left(3\frac{\nu_{G}^2}{\epsilon_f^2}+3M_G^2  \frac{\nu_{f}^2}{\epsilon_f^{2} }+ 
				3\left(\frac{M_{G}^2}{\epsilon_f^{4}}\right)\right)}.
			\tag{Setting B}
		\end{align}
		and the constants $\nu_f, \nu_{f,j}, \nu_G$, and $\nu_{G,j}$ are as specified in Lemma~\ref{cond-sec-mom}.
	\end{lemma}
	\subsection{Convergence Analysis} \label{sec:3.3}
	\begin{proposition}\label{bd-itercomp-sc-sgd}
		{Suppose $h$, defined as in \eqref{def-h}, is a convex function on an open set containing $\mathcal X$}. {Suppose Assumption ~\ref{dist-ass} holds and either \ieb{Assumption~\ref{dist-ass} or Assumption~\ref{dist-ass-3}} holds under Setting B. In addition, suppose Assumptions~\ref{sfo} and \ref{assump_sub} hold.} Consider the iterates
		generated by Algorithm~\ref{smooth scheme}. {If } $\bar{x}_{\widehat{K},K} \triangleq \tfrac{\sum_{k=\widehat{K}}^{K-1} \gamma_k x_k}{\sum_{k=\widehat{K}}^{K-1} \gamma_k},$ then for all $K>0$ and $\widehat{K}$ satisfying $0 \leq \widehat{K} < K-1$, 
		\begin{align}\label{bd-prop-h}
			\mathbb{E}\left[h(\bar{x}_{\widehat{K},K})-h(\x^*)\right] \leq
			\tfrac{\mathbb{E}[\|x_{\widehat{K}}-x^*\|^2]+\sum_{k={\widehat{K}}}^{K-1}\gamma_k^2 (M_G^2 + B^2)+\sum_{k={\widehat{K}}}^{K-1} \tfrac{\nu^2}{\sqrt{N_k}}}{\sum_{k={\widehat{K}}}^{K-1}2\gamma_k}.
		\end{align}
	\end{proposition}

	{We now present a rate statement for diminishing and constant steplengths.}
	\begin{theorem}[{\bf Rate statement for diminishing and constant steplengths}]\label{th-bd-itercomp-sc-sgd}
		{Suppose $h$, defined as in \eqref{def-h}, is a convex function on an open set containing $\mathcal X$}. {Suppose Assumption ~\ref{dist-ass} holds and either Assumption~\ref{dist-ass} or Assumption~\ref{dist-ass-3} holds under Setting B. In addition, suppose Assumptions~\ref{sfo} and \ref{assump_sub} hold.} Consider the iterates
		generated by Algorithm~\ref{smooth scheme}.
		
		\noindent (a) Suppose $\gamma_k = \tfrac{1}{k^{1/2+a}}$ and  $N_k \triangleq {\lceil {1/\gamma_k^4} \rceil}$ for all $k$ where $a < 1/2$. If $\widehat{K} \triangleq \lfloor K/2 \rfloor$, then the following holds for every integer $K \ge 2$.
		\begin{align}
			\mathbb{E}\left[(h(\bar{\x}_{\widehat{K},K})-h(\x^*))\right] & \leq
			(1/2-a)\tfrac{B^2 + \tfrac{1}{2a}(M_G^2 + B^2 + \nu^2)}{2(1-1/2^{1/2-a})} \tfrac{1}{K^{1/2-a}}.
		\end{align}
		\noindent (b) Given a positive integer $K$, suppose $\gamma_k \triangleq \sqrt{\tfrac{B^2}{(B^2+M_G^2+\nu^2)K}}$  and $N_k \triangleq {\lceil 1/\gamma_k^4 \rceil}$ for all $k$. Then the following holds. 
		\begin{align}
			\mathbb{E}\left[(h(\bar{\x}_{K})-h(\x^*))\right] & \leq \sqrt{\tfrac{(B^2+M_G^2+\nu^2)}{B^2 K}}.
		\end{align}
	\end{theorem}
	\begin{proof}
		\noindent (a) Suppose $\widehat{K} = \lfloor K/2\rfloor $ and $\gamma_k = \tfrac{\gamma_0}{k^{1/2+a}}$ for any $k \geq 0$. Then we have that
		\begin{align*}
			\sum_{k= \widehat{K}}^{K-1} \gamma_k & \geq \int_{\widehat{K}-1}^K \tfrac{1}{x^{1/2+a}} dx = \tfrac{K^{1/2-a}-(\widehat{K}-1)^{1/2-a}}{1/2-a} \geq  \tfrac{K^{1/2-a}- (K/2)^{1/2-a}}{1/2-a} \\ 
			\mbox{ and } \sum_{k= \widehat{K}}^{K-1} \gamma^2_k & \leq \int_{\widehat{K}-1}^{K} \tfrac{1}{x^{1+2a}} dx = \tfrac{K^{-2a}-(\widehat{K}-1)^{-2a}}{-2a} \leq \tfrac{(\widehat{K}-1)^{-2a}}{2a} \leq \tfrac{1}{2a} .
		\end{align*} 
		It follows that  
		\begin{align*}
			\mathbb{E}\left[(h(\bar{\x}_{\widehat{K},K})-h(\x^*))\right] & \leq
			\frac{B^2}{\sum_{k={\widehat{K}}}^{K-1}2\gamma_k } + \frac{\sum_{k={\widehat{K}}}^{K-1}\gamma_k^2 (M_G^2 + B^2)}{\sum_{k={\widehat{K}}}^{K-1}2\gamma_k }+ \frac{\sum_{k={\widehat{K}}}^{K-1} \tfrac{\nu^2}{\sqrt{N_k}}}{\sum_{k={\widehat{K}}}^{K-1}2\gamma_k} \\
			& \leq (1/2-a)\tfrac{B^2 + \tfrac{1}{2a}(M_G^2 + B^2 + \nu^2)}{2(1-(1/2)^{1/2-a})} \tfrac{1}{K^{1/2-a}}.
		\end{align*}
		
		\noindent (b) {Suppose $\widehat{K} = 0$} and $\gamma_k = \gamma$ for all $k$. Then we obtain the following bound. 
		\begin{align*}
			\mathbb{E}\left[(h(\bar{\x}_{K})-h(\x^*))\right] & \leq
			\tfrac{B^2}{2K\gamma} + \tfrac{(M_G^2+B^2+\nu^2) K \gamma^2 }{2 K \gamma} = 
			\tfrac{B^2}{2K\gamma} + \tfrac{(M_G^2+B^2+\nu^2)  \gamma }{2 }.
		\end{align*}
		By minimizing the right hand side, {which is a convex function} in $\gamma$, we obtain
		$$ -\tfrac{B^2}{2K\gamma^2} + \tfrac{(M_G^2+B^2+\nu^2)  }{2 } = 0 \implies \gamma^* = \sqrt{\tfrac{B^2}{(B^2+M_G^2+\nu^2)K}}.$$
		The resulting bound on the expected sub-optimality is
		\begin{align*}
			\mathbb{E}_{\tilde p}\left[(h(\bar{\x}_{K})-h(\x^*))\right] & \leq \sqrt{\tfrac{(B^2+M_G^2+\nu^2)}{B^2 K}}. \hspace{1in} \qed
		\end{align*}
	\end{proof}
	
	We now employ the aforementioned rate to compute the sample (or oracle) complexity of computing a random $\x_K$ such that $\mathbb{E}[h(\x_K)-h(\x^*)] \leq \epsilon$.   
	\begin{proposition}[{\bf Oracle complexity for diminishing \& constant steplengths}]\label{bd-itercomp-sc-sgd-2}
		{Suppose $h$, defined as in \eqref{def-h}, is a convex function on an open set containing $\mathcal X$}. {Suppose Assumption ~\ref{dist-ass} holds and either Assumption~\ref{dist-ass} or Assumption~\ref{dist-ass-3} holds under Setting B. In addition, suppose Assumptions~\ref{sfo} and \ref{assump_sub} hold.} Consider the iterates
		generated by Algorithm~\ref{smooth scheme}.
		
		\noindent (a) Suppose $\gamma_k = \tfrac{1}{k^{1/2+a}}$ and  $N_k \triangleq \lceil {1/\gamma_k^4} \rceil$ for all $k$ where $a < 1/2$. If $\widehat{K} \triangleq \lfloor K/2 \rfloor$, then the following holds for every integer $K \ge 2$. {Let $K({\epsilon})$ be any positive integer}, $K(\epsilon) \geq 2$, such that $\mathbb{E}_{\tilde p}[h(\x_{K(\epsilon)})-h(\x^*)] \leq \epsilon$. Then $\sum_{k=0}^{K(\epsilon)} N_k \leq \mathcal{O}(1/\epsilon^{(6+8a)/(1-a)}).$   
		
		\noindent (b) Given a positive integer $K$, suppose $\gamma_k \triangleq \sqrt{\tfrac{B^2}{(B^2+M_G^2+\nu^2)K}}$  and $N_k \triangleq \lceil {1/\gamma_k^4} \rceil$ for all $k$. 
		Let $K_{\epsilon}$ be any positive integer $K(\epsilon) \geq 2$ such that $\mathbb{E}_{\tilde p}[h(\x_{K(\epsilon)})-h(\x^*)] \leq \epsilon$. Then $\sum_{k=0}^{K(\epsilon)} N_k \leq \mathcal{O}(1/\epsilon^6).$   
	\end{proposition}
	\begin{proof}
		\noindent (a). By utilizing Theorem~\ref{th-bd-itercomp-sc-sgd}(a), we have that 
		$K(\epsilon) = \lceil \tfrac{\widehat{D}}{\epsilon^{2/(1-a)}}\rceil$, {where $\hat D>0$}. Consequently, we have that 
		\begin{align*}
			\sum_{k=0}^{K(\epsilon)} N_k & \leq \sum_{k=0}^{K(\epsilon)} ({k+1})^{2+4a} = \sum_{t=1}^{K(\epsilon)+1} t^{2+4a} \leq \int_1^{K(\epsilon)+2} x^{2+4a} dx\\
			& \leq \frac{(K(\epsilon)+2)^{3+4a}}{3+4a} \leq \frac{{2^{3+4a}}\widehat{D}^{3+4a}}{(3+4a)\epsilon^{(6+8a)/(1-a)}} .
		\end{align*}
		
		\noindent (b) By utilizing Theorem~\ref{th-bd-itercomp-sc-sgd}(b), we have that 
		$K(\epsilon) = \lceil \tfrac{\widehat{D}}{\epsilon^{2}}\rceil$, {implying that $\gamma_k = \gamma = \tfrac{\tilde{c}}{\sqrt{K}}$ where $\tilde{c} \triangleq \sqrt{\tfrac{B^2}{(B^2+M_G^2+\nu^2)}}$. It follows that $N_k = N = \lceil \tfrac{K_{\epsilon}^2}{\tilde{c}^4}\rceil.$ The oracle complexity may then be bounded as }
		$ \sum_{k=0}^{K(\epsilon)} N_k \leq \frac{8\widehat{D}^{3}}{(3\tilde{c}^4)\epsilon^{6}}.$
	\end{proof}
	
	Next, we prove a.s. convergence of the sequence to a solution $\x^*$. 
	\begin{theorem}[Almost sure convergence]\label{th-asconv-itercomp-sc-sgd}
		{Suppose $h$, defined as in \eqref{def-h}, is a convex function on an open set containing $\mathcal X$}. {Suppose Assumption ~\ref{dist-ass} holds and either Assumption~\ref{dist-ass} or Assumption~\ref{dist-ass-3} holds under Setting B. In addition, suppose Assumptions~\ref{sfo} and \ref{assump_sub} hold.} Consider the iterates
		generated by Algorithm~\ref{smooth scheme}, where $\sum_{k=0}^{\infty} \gamma_k=\infty$, $\sum_{k=0}^{\infty} \gamma_k^2 < \infty$, and $\sum_{k=0}^{\infty} \tfrac{1}{\sqrt{N_k}} < \infty.$ Then $\x_k \xrightarrow[a.s.]{k \to \infty} \mathcal X^*.$ 
	\end{theorem}
	\begin{proof} We resume our argument utilizing the following inequality.  
		\begin{align*}
			{1\over 2}\|\x_{k+1}-\x^*\|^2& \leq {1\over 2}\|\x_k-\gamma_k \us{(d_k+\bar{w}_k)}-x^*\|^2\\
			&={1\over 2}\|\x_k-\x^*\|^2+{1\over 2}\gamma_k^2\|\us{d_k + \bar{w}_k}\|^2-\gamma_k(x_k-x^*)^T(\us{d_k}+\bar w_{k}) \\
			&\leq {1\over 2}\|\x_k-\x^*\|^2+{1\over 2}\gamma_k^2\|\us{d_k + \bar{w}_k}\|^2-\gamma_k(x_k-x^*)^T(\us{d_k}) \\
			& + \gamma_k^2 \|x_k-x^*\|^2 +\|\bar w_{k}\|^2. 
		\end{align*}
		This implies the following inequality holds. 
		\begin{align*}
			\|\x_{k+1}-\x^*\|^2 & \leq \|\x_k-\x^*\|^2+\gamma_k^2\|\us{d_k + \bar{w}_k}\|^2-2\gamma_k (h(\x_k)-h(\x^*)) \\ &+ \gamma_k^2 \|x_k-x^*\|^2 +\|\bar w_{k}\|^2 \\
			& = (1+\gamma_k^2) \|\x_k-\x^*\|^2+\gamma_k^2\|\us{d_k + \bar{w}_k}\|^2\\ &-2\gamma_k (h(\x_k)-h(\x^*)) +\|\bar w_{k}\|^2. 
		\end{align*}
		By taking expectations conditioned on $\mathcal{F}_k$, we have the following inequality. 
		\begin{align*}
			\mathbb{E}[\|\x_{k+1}-\x^*\|^2 \mid \mathcal F_k] & \leq (1+\gamma_k^2) \|\x_k-\x^*\|^2+\gamma_k^2M_G^2-2\gamma_k (h(\x_k)-h(\x^*)) + \tfrac{\nu^2}{\sqrt{N_k}}. 
		\end{align*}
		Since $\sum_{k=0}^{\infty} \gamma_k^2 < \infty$, $\sum_{k=0}^{\infty} \tfrac{1}{\sqrt{N_k}} < \infty$, it follows that $\{\|\x_k-\x^*\|^2\}$ is a convergent sequence in an a.s. sense and $\sum_{k=0}^{\infty}\gamma_k (h(\x_k)-h(\x^*)) < \infty$ a.s. Consequently, $\{\x_k\}$ is bounded a.s. and has a convergent subsequence, indexed by $\mathcal{I}$. Since $\sum_{k=0}^{\infty}\gamma_k (h(\x_k)-h(\x^*)) < \infty$ a.s. and $\sum_{k=0}^{\infty} \gamma_k = \infty$, it follows that $\liminf_{k \to \infty, k \in \mathcal I} h(\x_k) = h(\x^*)$ a.s. Consequently, there exists a subsequence of $\{\x_k\}$ that converges to the solution set $\mathcal X^*$ almost surely. But we have that $\{\|\x_k-\x^*\|^2\}$ is convergent a.s. and converges to zero along some subsequence. Consequently, the entire sequence $\{\|\x_k-\x^*\|^2\}$ converges to zero a.s. and the result holds. 
	\end{proof}
	\section{Numerical Results}\label{sec:4}
	{In this section, we compare the performance of our scheme with standard stochastic approximation and integer programming approaches on two sets of examples.  In all instances,  the components of $ x_0 \in \mathbb{R}^n$ are chosen randomly from the standard uniform distribution. In standard (\textbf{SA}) and (\textbf{batch-SA}) algorithms, the step length sequence is $\{ 1/\sqrt{k}\} $, while in (\textbf{batch-SA}), we compute the approximate subgradients using batch size of $ 100 $ samples.}
	{In {\bf (r-VRSA)}, the step length sequence is $ \{ \tfrac{\gamma_0}{k^{1/2+a}}\} $. The parameters $ \gamma_0 $ and $ a $ \uss{are} chosen as $ \gamma_0=20,\ a=0 $ in Example 1 and $ \gamma_0=1,\ a=0 $ in Example 2. \uss{In all instances,  the components of $ x_0
			\in \mathbb{R}^n$ are chosen randomly from the standard uniform distribution.}} 
	
	\noindent \textbf{Example 1. Set Covering. }  {Demand} is assumed to be uniformly distributed on \ieb{$\Kscr =\{\zeta \mid \|\zeta-\alpha\| \leq \alpha\}$} while the cost of operating a vehicle on route $j$ is given by $c_j$ while {$ \beta $ is a cost threshold.} By Prop \ref{bound-asym-sett_B}, we may rewrite the problem (\ref{ex1}) as
	\begin{align*}
		\min_{\x} & \  h(\x) \ \triangleq \ -\log\mathbb{E}_{\tilde p(\xi)} [F(\x,\xi)],  \quad \afj{\mbox{s. t. }} \ c^T \x \leq {\beta}, \quad\ \x \geq 0 \mbox{ where } 
	\end{align*}
	$F(\x,\xi)  \triangleq  \mathcal{C}(2\pi\ieb{\sigma^2})^{\ieb{d}/2} 
	e^{-g(\ieb{\x},\xi)+\frac{\|\xi\|^2_{2}}{2\sigma^2}}$, $ \mathcal{C}  \triangleq \tfrac{1}{\mathrm{Vol}(\Kscr)} \ieb{\tfrac{1}{\Gamma(1+d/2)}}$, $ \tilde p(\xi) \triangleq 
	\tfrac{1}{ (2\pi\ieb{\sigma^2})^{\ieb{d}/2}  } e^{\frac{-\|\xi\|^2_{2}}{2\sigma^2}},$
	\mbox{and,}
	$g(\ieb{\x},\xi) \triangleq \max\left\{ \frac{1}{\alpha^2}\|\xi\|_{p}^2,\left(\frac{\max(\xi_1,0)}{T_{1,\bullet}\x-\mu_1}\right)^2,\cdots, \left(\frac{\max(\xi_d,0)}{T_{d,\bullet}\x-\mu_d}\right)^2  \right\}.$
	
	Note that $ \xi's $ are normally distributed with zero mean and standard deviation $ \sigma $ where  $ \sigma^2 = \alpha^2  $. We compare
	the performance of these algorithms for different setups.  In these problems,
	the vehicle routing network is randomly generated and corresponding incidence
	matrix $ T \in \mathbb{R}^{d\times n} $ is obtained. The elements of cost
	vector $ c \in \mathbb{R}^n$ are randomly chosen from the uniform distribution
	on $ [0,c_{\max}] $. \uss{We compare the performance and quality of
		the solutions with those obtained via an integer programming approximation as
		proposed in~\cite{luedtke08sample}. This avenue employs a sample average
		approximation approach facilitated by integer programming (denoted by
		(SAA-IP)), defined as follows: 
		\begin{align} \tag{SAA-IP$_N$} \begin{aligned}
				\max_{\x \in \mathcal{X}, \mathbf{z} \in \{0,1\}^N} & \   \frac{1}{N} \sum_{j=1}^N z_j \\
				\mbox{subject to }  & \, T_{i,\bullet}\x \geq v_i, && i = 1, \hdots, d\\
				&  v_i \geq \zeta_{i}^{j}z_j, && i = 1, \hdots, d, \ j = 1, \hdots, N\\
				& v_i \geq 0. && i = 1, \hdots, d
			\end{aligned} 
		\end{align}
		In this formulation, auxiliary variables $ v_i $ for $ i = 1, \hdots, d $
		represent $ T_{i,\bullet}\x $ and $ z_j = 1$ (or $0$), then the constraints $
		T_{i,\bullet}\x \geq \zeta_{i}^{j} $ for $ i = 1, \hdots, d $ corresponding to
		the realization $ j $ in the sample are enforced (or not enforced). We solve
		this problem using Gurobi MIP solver. The sample size for SAA-IP scheme is $ N
		= 1e4 $.}
	
	To  compare across the solutions of various schemes including  SA, batch-SA, r-VRSA, and SAA-IP, we first generate samples of demand vector $ \zeta $ from the set \ieb{$\Kscr =\{\zeta \mid \|\zeta-\alpha\| \leq \alpha\}$}. We then use Monte Carlo simulation to estimate $ f(\x) \ \triangleq \ \mathbb{P}\{\zeta \in \Kscr \mid T \x \geq \zeta\} $ for the solution $ \x $ of each scheme.  In Table \ref{tab:1}, the
	first column prescribes problem parameters as follows: (Problem \#, $d$, $n$,
	$\alpha $, $ c_{\max} $, $ \gamma $).
	\ieb{In SAA-IP scheme, the \textit{Gap} refers to the reported gap between upper and lower bounds.} Note that the table shows the probability being maximized.

\begin{table}[htbp]
	\resizebox{\textwidth}{!}{%
		\begin{tabular}{|c|cc|cc|cc|cc|}
			\hline
			\multirow{4}{*}{\textbf{Problem}} & \multicolumn{6}{c|}{\multirow{2}{*}{\textbf{r-VRSA}}}                                                           & \multicolumn{2}{c|}{\multirow{2}{*}{\textbf{SAA-IP}}} \\
			& \multicolumn{6}{c|}{}                                                                                           & \multicolumn{2}{c|}{}                                 \\ \cline{2-9} 
			& \multicolumn{2}{c|}{\textbf{B=1e6}} & \multicolumn{2}{c|}{\textbf{B=1e7}} & \multicolumn{2}{c|}{\textbf{B=1e8}} & \multicolumn{2}{c|}{\textbf{B=1e4}}                   \\ \cline{2-9} 
			& \textbf{f(x)}    & \textbf{Time}    & \textbf{f(x)}    & \textbf{Time}    & \textbf{f(x)}    & \textbf{Time}    & \textbf{f(x)}              & \textbf{Gap}             \\ \hline
			(1, 10, 9, 10, 5, 170)            & 0.9840           & 24s              & 0.9847           & 154s             & 0.9852           & 821s             & 0.9852                     & \%0                      \\
			(2, 14, 16, 8, 3, 46)             & 0.8341           & 26s              & 0.8356           & 158s             & 0.8357           & 1248s            & 0.8346                     & \%0.4                  \\
			(3, 18, 23, 16, 7, 250)                & 0.9325           & 32s              & 0.9327           & 172s             & 0.9328           & 1335s            & 0.9317                     & \%0.2                    \\
			(4, 23, 54, 40, 20, 530)               & 0.8255           & 33s              & 0.8767           & 177s             & 0.8768           & 1391s            & 0.8759                     & \%0.8                  \\ \hline
	\end{tabular}}
	\caption{Set covering problem. (SAA-IP algorithm is terminated after 10000s.)}
	\label{tab:1}

\end{table}
\noindent \textbf{Example 2. \ieb{Robust portfolio selection problem.}} {We now consider the \ieb{robust} portfolio selection problem described in Section~\ref{sec:app}. We compare the proposed approach with the quadratic minimization (QM) framework~\cite{bardakci2019distributionally} through which exact solutions are available}. The portfolio weights are restricted to lie in the set  $ \mathcal{X}$ where $\mathcal{X} \triangleq \left\{\x: {\bf 1}^T \x =1 \text{ and } \x \geq 0 \right\}$. The parameter $ \alpha $ is set as $( \alpha = 0 )$.  Given a threshold $\alpha$ and an allocation $\x$, we use the proposed framework to estimate the probability of a loss being less or equal than $\alpha$ as
$f_{\alpha}(\x) \triangleq \mathbb{P} \left\{\tilde{\zeta}: \tilde{\zeta}^T\x\leq -\alpha\right\}$ where $ \tilde{\zeta} = \zeta + \boldsymbol{\mu} $.  In our simulations, given the number of assets $ n $, mean $ \mathbf{\mu} $ (randomly generated), and covariance of random returns $ \Sigma $,  $ \zeta $, the returns, are assumed to be uniformly distributed over the set
$\Kscr_{\epsilon} = \{ \zeta \in \mathbb{R}^n:  \zeta^T \Sigma^{-1}  
\zeta \leq 1\}.$
In Table \ref{tab:port_table}, Problem column corresponds to (Problem no., number of assets $ n $, $ \gamma $ ).

\begin{figure}[htbp]
	\vspace{-0.1in}
	\centering
	\subfigure[Example 1 (Vehicle Routing)]{\includegraphics[width=8cm]{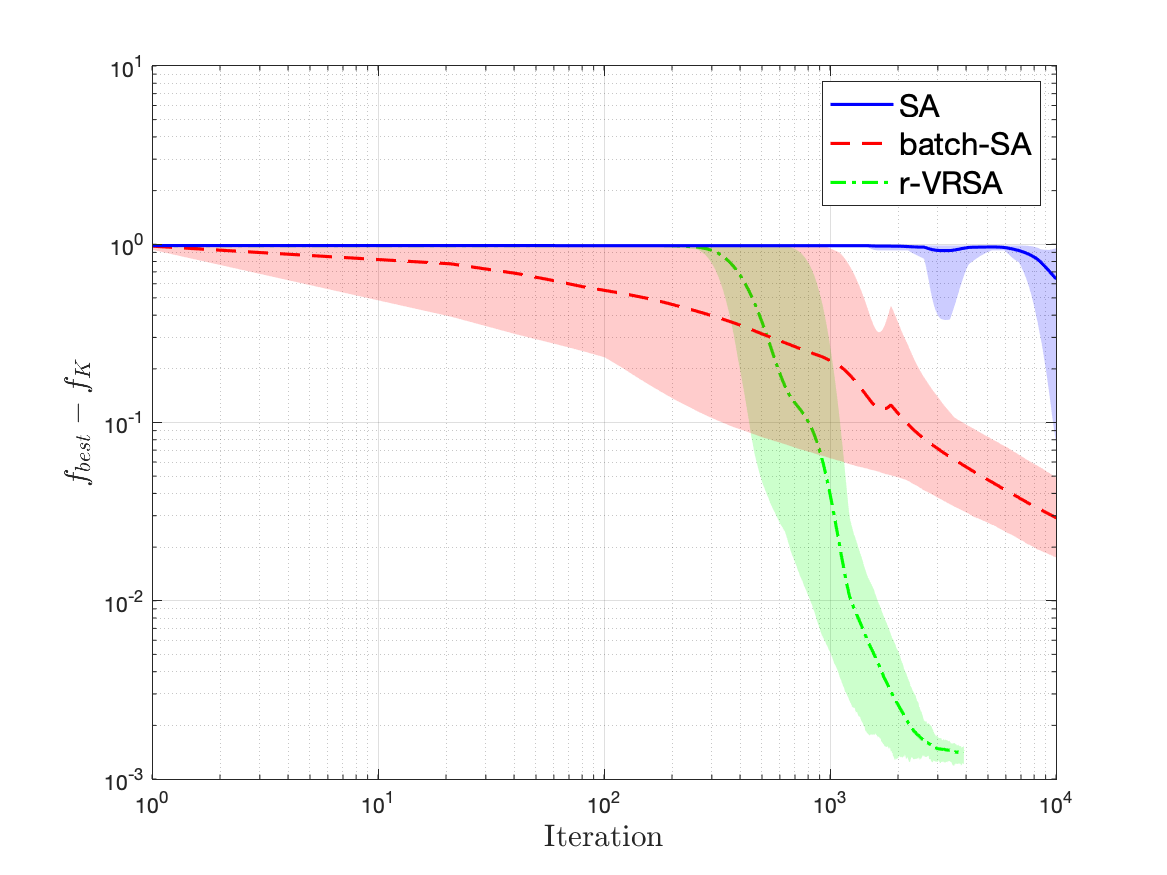}\label{fig:vrt}}
	\subfigure[Example 2 (Portfolio Optimization)]{\includegraphics[width=8cm]{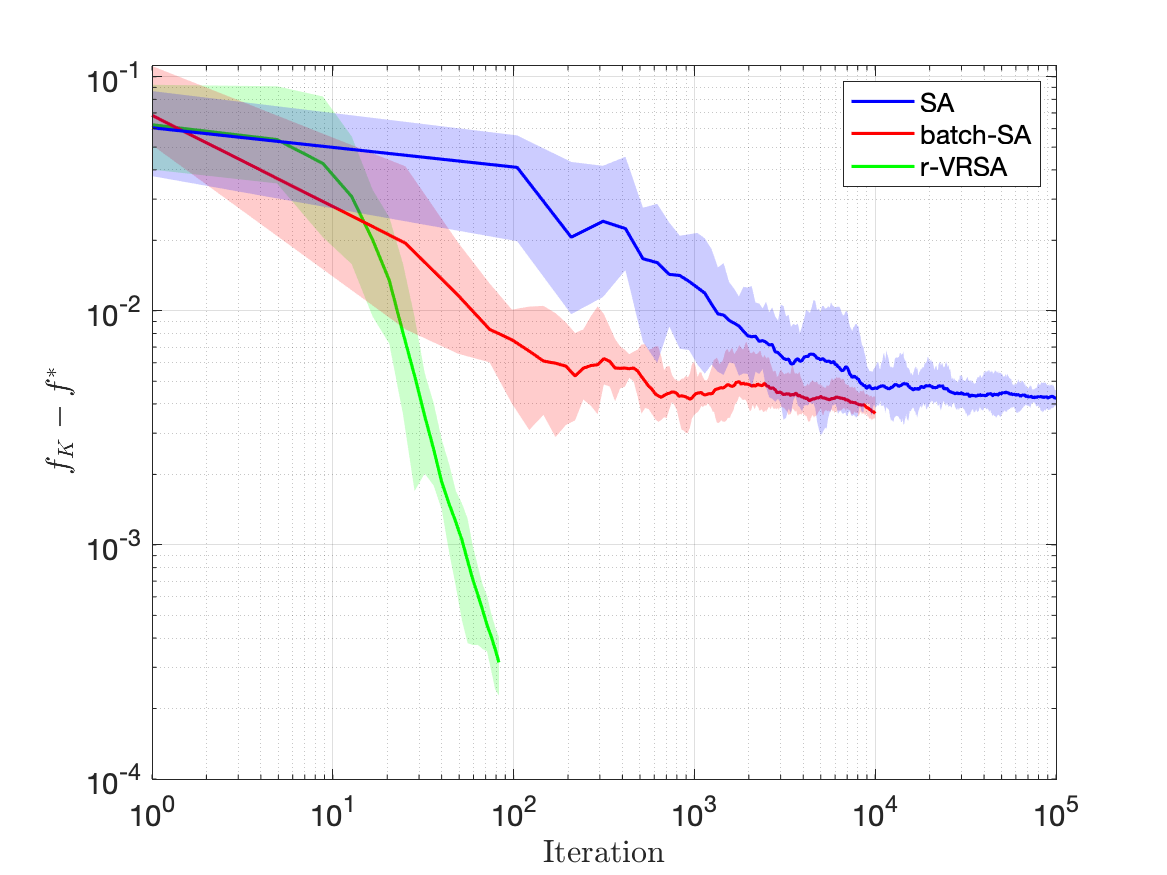}\label{fig:port}}
	\vspace{-0.2in}
	\caption{Comparison of algorithms.}\label{fig:trajectory}

\end{figure}

{In Figure~\ref{fig:vrt}, the budget is $1e7$ and $ f_{best} $ (an approx. of $f^*$) is obtained by running the {\bf (r-VRSA)} with a budget of $1e11$. In Figure~\ref{fig:port}, the budget is $1e6$. In both figures, the standard (SA) algorithm is terminated after $1e5$ iterations.}

\begin{table}[htbp]
	\centering
	\small
	\begin{tabular}{|c|cc|cccc|c|}
		\hline
		\multirow{2}{*}{\textbf{Problem}} & \multicolumn{2}{c|}{\textbf{batch-SA}} & \multicolumn{4}{c|}{\textbf{r-VRSA}}                              & \multirow{2}{*}{\textbf{\begin{tabular}[c]{@{}c@{}}QM\\ ($f^{*}$)\end{tabular}}} \\ \cline{2-7}
		& \textbf{B=1e5}     & \textbf{B=1e6}    & \textbf{B=1e5} & \textbf{B=1e6} & \textbf{B=1e7} & \textbf{B=1e8} &                                                                                  \\ \hline
		(1, 4, 0.25)                      & 0.3730             & 0.3723            & 0.3715         & 0.3712         & 0.3712         & 0.3711         & 0.3710                                                                           \\
		(2, 16, 0.2)                      & 0.3000             & 0.2993            & 0.2976         & 0.2964         & 0.2963         & 0.2961         & 0.2961                                                                           \\
		(3, 64, 0.05)                     & 0.3818             & 0.3752            & 0.3894         & 0.3787         & 0.3751         & 0.3743         & 0.3743                                                                           \\
		(4, 128, 0.15)                    & 0.0899             & 0.0872            & 0.0992         & 0.0886         & 0.0869         & 0.0868         & 0.0867                                                                           \\
		(5, 256, 0.1)                     & 0.1321             & 0.0991            & 0.1340         & 0.1031         & 0.0975         & 0.0972         & 0.0966                                                                           \\ \hline
	\end{tabular}
	\caption{Portfolio selection}
	\label{tab:port_table}
\end{table}

\noindent {\bf Comments.} {Several observations can be made. (i) In Example 1,
	{\bf (r-VRSA)} obtains near-optimal solutions within $1$-$2\%$ of the  time taken
	by {\bf (SAA-IP)}, an integer programming approach. (ii) While {\bf (batch-SA)}
	performs reasonably in Setting A, it tends to degenerate in Setting B. Further, convergence theory is unavailable for
	this scheme. Such schemes perform less favorably in comparison to {\bf
		(r-VRSA)}.  (iii) ({\bf SAA-IP}) produces solutions of inferior gap as
	dimension grows and cannot accommodate growing number of samples, impacting
	solution quality.}

\vspace{-0.2in}
\section{Concluding remarks and future work} {Traditional approaches for contending with
	chance-constrained optimimzation problems have relied on resolving convex
	approximations or computing stationary points. We concentrate our efforts on a
	subclass of such problems that require maximization the probability of a
	suitably specified event. By leveraging a recent result on non-Gaussian
	integrals of PHFs, we show that the probability of interest is an expectation
	of a possibly nonsmooth integrand. It is then shown that the composition of
	this expectation with a suitably specified smooth convex function leads to a
	convex program. However, a direct application of SA schemes is impeded by the
	inability to generate unbiased samples of the gradient. This motivated the
	development of a regularized variance-reduced SA scheme {\bf (r-VRSA)} that matches the
	optimal rate of subgradient methods for nonsmooth convex optimization problems
	but has somewhat poorer sample complexity, a consequence of the unavailability of conditionally unbiased gradients.  We believe that this set of contributions represents amongst the first avenues (to the best of our knowedge) for tractably resolving probability maximization problems and is a crucial
	first step in examining  more intricate problems in 
	chance-constrained optimization.  } 

This framework will provide the cornerstone for at least two key generalizations in our future work. 
\begin{enumerate}
	\item[(i)] {\em Extensions to log-concave measures.} First, this avenue may be extended to symmetric log-concave measures,  subsuming Guassian, Laplace, Subbotin, amongst others. 
	\cl{Pathways being exploited include alternative representations of probability density functions such as the so-called layer cake representation; e.g., see~\cite{lieb2001analysis}.}
	
	\item[(ii)] {\em Extensions to chance-constrained regimes.} This framework also allows for contending with constrained regimes. Consider the following chance-constrained problem and its expectation-valued counterpart.
	\begin{align*}
		\left\{\begin{aligned}
			\min_{\x \in \Xscr} & \quad  f(\x) \\
			\mbox{subject to } & \quad \mathbb{P}\{ \zeta \mid \zeta \in \K(\x)\} \geq \bar{p}
		\end{aligned}
		\right\} \equiv
		\left\{\begin{aligned}
			\min_{\x \in \Xscr} & \quad  f(\x) \\
			\mbox{subject to } & \quad \psi(\mathbb{E}[F(\x,\xi)]) \leq \bar{c},
		\end{aligned}
		\right\}
	\end{align*}
	where $\bar{c}$ is related to $\bar{p}$.  Assuming that $f$ is a convex function on $\Xscr$, we observe that the techniques in this paper allow for recasting the chance
	constrained problem can be recast as a convex optimization problem with
	nonsmooth compositional expectation-valued constraints.
\end{enumerate}
Extensions and generalizations captured in (i) and (ii), while challenging, remain the focus of future work.\\
}

\vspace{-0.2in}
\section{Appendix}
\textbf{Proof of Theorem~\ref{theo-convexity-settingA}:}
(a) \uss{When considering uniform distributions over a compact and convex  set
$\mathcal{K}$, the density is constant in this set and zero outside the set. It
can then be concluded that $\zeta$  has a log-concave density. Furthermore,
$\zeta$ has a symmetric density about the origin since $\mathcal{K}$ is a
symmetric set about the origin. Hence by Lemma 6.2 in \cite{bobkov2010}, $h$ is
convex where $h(\x)\triangleq 1/f(\x)$.}\\

\noindent (b) Since \eqref{hx} is a convex program, any solution $\x^*$ satisfies
${ h(\x^*) \leq h(\x), \ \forall \mathbf{x} \in \mathcal{X}.}$
From the positivity of $f$ over $\mathcal X$, 
$ \frac{1}{f(\x^*)} \leq \frac{1}{f(\x)}$ for every $\x \in \mathcal{X}$ implying that $f(\x^*) \geq f(\x)$ for every $\x \in \mathcal{X}.$
Consequently, $\x^*$ is a global maximizer of \eqref{hx}.\\ \qed 

\noindent \textbf{Proof of Lemma~\ref{unimodal}:}
\us{We prove this result by showing the unimodality of $f$ on $\Real_+$ where $f(u) = u^c e^{-u}$, implying that $f'(u) = cu^{c-1}e^{-u} - u^ce^{-u} = 0$ if $u = c.$ Furthermore, $f'(u) > 0$ when $u < c$ and $f'(u) < 0$ when $u > c$. Finally, $f(0) = 0$. It follows that $u^* = c$ is a maximizer of $u^ce^{-u}$ on $[0,\infty)$ where $f(c) = \frac{c^c}{e^{c}}$.}\\  \qed

\noindent \textbf{Proof of Proposition~\ref{bound sub}:}
\uss{Recall the definition of $F(\x,\xi)$ from the statement of Lemma~\ref{lemma_Clarke_F}.}
We prove (a) by considering two cases. Case (i): $\xi \in \Xi_1(\x) \cup \Xi_0(\x).$ It follows that 
\begin{align*} |F(\x,\xi)|^2 = \mathcal{C}^2_{\Kscr} 
\left( (2\pi)^{{n}} e^{-2\mid\xi^T\x \mid^2+\|\xi\|_{\Kscr}^2}  \right) 
& \leq \mathcal{C}^2_{\Kscr} \left( (2\pi)^{{n}} e^{-2\mid\xi^T\x \mid^2+\mid\xi^T\x \mid^2}  \right)  \leq \mathcal{C}^2_{\Kscr} (2\pi)^n.
\end{align*}
Case (ii): $\xi \in \Xi_2(\x).$ Proceeding similarly, we obtain that
\begin{align*}
|F(\x,\xi)|^2
& \leq \mathcal{C}^2_{\Kscr} \left( (2 \pi)^{{n}} e^{\us{-2\|\xi\|_{\Kscr}^2}+{\us{\|\xi\|^2_{\Kscr}}} } \right) \leq  
\mathcal{C}^2_{\Kscr} \left( (2 \pi)^{{n}} e^{\us{-2\|\xi\|_{\Kscr}^2}+{\us{\|\xi\|^2_{\Kscr}}} } \right) \leq \mathcal{C}^2_{\Kscr} (2\pi)^n.
\end{align*}
Consequently, $|F(\x,\xi)|^2 \leq \mathcal{C}^2_{\Kscr} (2\pi)^n$ for every $\xi \in\Real^n$.  

\noindent (b) \uss{We observe that $\partial F(\x,\xi)$ is defined as follows.
\begin{align*} \scriptsize
	\partial F(\x,\xi) = \begin{cases}  \left( \ieb{\mathcal{C}_{\Kscr}}(2\pi)^{n/2} (-2\xi \xi^T\x) e^{-|\xi^{\intercal} \x |^{\us{2}}+\frac{\us{\|\xi\|_{\Kscr}^2}}{2}}\right),  & \xi \in \Xi_1(\x) \triangleq \left\{\xi \mid |\xi^{\intercal} \x |^{\us{2}} > \|\xi\|^{2}_{\mathcal{K}}\right\} \\ 
		\left( -\ieb{\mathcal{C}_{\Kscr}} (2\pi)^{n/2}  e^{-\max\{|\xi^{\intercal} \x |^{\us{2}}, {\|\xi\|^{\us{2}}_{\mathcal{K}}}\}+\frac{\us{\|\xi\|_{\Kscr}^2}}{2}}\right) \left[0,2\xi (\xi^T\x)\right],  & \xi \in \Xi_0(\x) \triangleq \left\{\xi \mid |\xi^{\intercal} \x |^{\us{2}} = \|\xi\|^{2}_{\mathcal{K}}\right\}\\
		{\bf 0}.  & \xi \in \Xi_2(\x) \triangleq \left\{\xi \mid |\xi^{\intercal} \x |^{\us{2}} < \|\xi\|^{2}_{\mathcal{K}}\right\}
	\end{cases}
\end{align*}
Consequently, it follows that $\mathbb{E}_{\ieb{\tilde{p}}} \left[\|G(\x,\xi)\|^2\right]$ is bounded as follows.
\begin{align}\notag
	& \quad \mathbb{E}\left[\|G(\x,\xi)\|^2\right]  = \int_{\Xi} \|G(\x,\xi)\|^2 \tilde{p}(\xi) d\xi \\
	& = \int_{\Xi_1(\x)} \|G(\x,\xi)\|^2 \tilde{p}(\xi) d\xi \notag 	
	+ \int_{\Xi_2(\x)} \|\underbrace{G(\x,\xi)}_{ \ = \ {\bf 0}}\|^2 \tilde{p}(\xi) d\xi 	
	\\ & + \int_{\Xi_0(\x)} \|G(\x,\xi)\|^2 \tilde{p}(\xi) d\xi \label{decomp-int}\\
	\notag			&  = \int_{\Xi_1(\x)} \|G(\x,\xi)\|^2 \tilde{p}(\xi) d\xi, 	
\end{align} 
where the last equality follows from observing that $G(\x,\xi) = 0$ for $\xi \in \Xi_2(\x)$ and the integral in \eqref{decomp-int} is zero because $\Xi_0(\x)$ is a measure zero set. 
It follows that $\mathbb{E}[\|G(\x,\xi)\|^2]$ can be bounded as follows:
\begin{align} 
	& \quad \mathbb{E}[\|G(\x,\xi)\|^2]\notag \\
	& = \int_{\Xi_1(\x)}\left(\mathcal{C}^2_{\Kscr}(2\pi)^{n} (4\|\xi\|_2^2 (\xi^T\x)^2) e^{-2(\xi^T\x)^2 +\us{\|\xi\|^2_{\Kscr}}} \right)\frac{1}{D_{\Kscr}}(2\pi)^{-\frac{n}{2}} e^{\frac{\us{\|\xi\|^2_{\Kscr}}}{2}} d\xi \label{eq:bound} \\
	& \leq \int_{\Xi_1(\x)}\left( \mathcal{C}^2_{\Kscr}(2\pi)^{n} (4\|\xi\|_2^2 (\xi^T\x)^2) e^{-(\xi^T\x)^2}  \right)\frac{1}{D_{\Kscr}}(2\pi)^{-\frac{n}{2}} e^{\frac{\us{\|\xi\|^2_{\Kscr}}}{2}} d\xi, 
\end{align}
where the inequality follows from $\xi \in \Xi_1(\x,u)$. Next, we consider the expression $(\xi^T\x)^2 e^{-(\xi^T\x)^2}$ or $ue^{-u}$. We note that by Lemma~\ref{unimodal}, $ue^{-u}$ is a unimodal function and $u^* = 1$ is a maximizer with value $e^{-1}$. \us{Consequently, we have that  $$\max_{\{(\xi^T\x) \mid \xi \in \Xi(\x)\}} (\xi^T\x)^2 e^{-(\xi^T\x)^2} \leq \max_{u \in \Real_+} u e^{-u}\overset{\tiny \mathrm{lemma}~\ref{unimodal}}{\leq} \tfrac{1}{e}, $$ implying}  that 
\begin{align*}
	& \quad	\mathbb{E}[\|G(\x,\xi)\|^2]   \leq \int_{\Xi_1(\x)}\left(\mathcal{C}_{\Kscr}^2(2\pi)^{n} (\|\xi\|_2^2 (\xi^T\x)^2) e^{-(\xi^T\x)^2} \right)\tfrac{1}{D_{\Kscr}}(2\pi)^{-\tfrac{n}{2}}e^{-\tfrac{\us{\|\xi\|^2_{\Kscr}}}{2}} d\xi \\
	& \leq  e^{-1} \mathcal{C}^2_{\Kscr}(2\pi)^{n}\int_{\Xi_1(\x)} \|\xi\|_2^2  \tfrac{1}{D_{\Kscr}}(2\pi)^{-\tfrac{n}{2}} e^{-\tfrac{\us{\|\xi\|^2_{\Kscr}}}{2}} d\xi \\
	& \leq e^{-1} \mathcal{C}^2_{\Kscr}(2\pi)^{n}\int_{\Real^n} \|\xi\|_2^2   \tfrac{1}{D_{\Kscr}}(2\pi)^{-\tfrac{n}{2}} e^{-\tfrac{\us{\|\xi\|^2_{\Kscr}}}{2}} d\xi 
	=  e^{-1} \mathcal{C}^2_{\Kscr}(2\pi)^{n}\mathbb{E}_{\tilde p}[\|\xi\|_2^2]. \hspace{0.51in}  
\end{align*}
}
\qed

\noindent \textbf{Proof of Proposition~\ref{bound-lp-ball-sett_A}:}
\ieb{(a) Since $\|\xi\|_{\Kscr}^2 = \|\xi\|_p^2$, it follows from Theorem~\ref{theo-rep} that 
\begin{align*}
	f(\x)&  =  \int_{\mathbb{R}^n} \left(\mathcal{C}(2\pi\sigma^2)^{\tfrac{n}{2}} e^{-\max \{ \mid \xi^T \x \mid^2, \us{\|\xi\|_{p}^2}\} } \right) d\xi\\
	& = \int_{\mathbb{R}^n} \underbrace{\left(\mathcal{C}(2\pi\sigma^2)^{\tfrac{n}{2}} e^{-\max \{ \mid \xi^T \x \mid^2, \us{\|\xi\|_{p}^2} \}+\tfrac{\us{\|\xi\|_{2}^2}}{2\sigma^2} } \right)}_{\triangleq F(\x,\xi)} \underbrace{(2\pi\sigma^2)^{-\tfrac{n}{2}} e^{-\tfrac{\us{\|\xi\|_{2}^2}}{2\sigma^2}}}_{\triangleq \tilde{p}(\xi)} d\xi,
\end{align*}

\noindent (b) Omitted (similar to proof of Proposition~\ref{bound sub}(a).

\noindent (c) Next, we derive a bound on the second moment of $\|G(\x,\xi)\|$ akin to Prop.~\ref{bound sub}(b). We observe that $ \partial F(\x,\xi) $ is defined as
\begin{align*} \scriptsize
	\partial F(\x,\xi) = \begin{cases}  \left( \mathcal{C}(2\pi\sigma^2)^{n/2} (-2\xi \xi^T\x) e^{-|\xi^{\intercal} \x |^{\us{2}}+\frac{\us{\|\xi\|_{2}^2}}{2\sigma^2}}\right),  & \xi \in \Xi_1(\x) \triangleq \left\{\xi \mid |\xi^{\intercal} \x |^{\us{2}} > \|\xi\|^{2}_{p}\right\} \\ 
		\left( -\mathcal{C} (2\pi\sigma^2)^{n/2}  e^{-\max\{|\xi^{\intercal} \x |^{\us{2}}, {\|\xi\|^{\us{2}}_{p}}\}+\frac{\us{\|\xi\|_{2}^2}}{2\sigma^2}}\right) \left[0,2\xi (\xi^T\x)\right],  & \xi \in \Xi_0(\x) \triangleq \left\{\xi \mid |\xi^{\intercal} \x |^{\us{2}} = \|\xi\|^{2}_{p}\right\}\\
		{\bf 0}.  & \xi \in \Xi_2(\x) \triangleq \left\{\xi \mid |\xi^{\intercal} \x |^{\us{2}} < \|\xi\|^{2}_{p}\right\}
	\end{cases}
\end{align*}
Consequently, $\mathbb{E}[\|G(\x,\xi)\|^2]$ can be bounded as follows. 
\begin{align}\notag
	& \quad \mathbb{E}\left[\|G(\x,\xi)\|^2\right]  = \int_{\Xi} \|G(\x,\xi)\|^2 \tilde{p}(\xi) d\xi \\
	& = \int_{\Xi_1(\x)} \|G(\x,\xi)\|^2 \tilde{p}(\xi) d\xi \notag 	
	+ \int_{\Xi_2(\x)} \|\underbrace{G(\x,\xi)}_{ \ = \ {\bf 0}}\|^2 \tilde{p}(\xi) d\xi 	
	\\ & + \int_{\Xi_0(\x)} \|G(\x,\xi)\|^2 \tilde{p}(\xi) d\xi \label{decomp-int2}\\
	\notag			&  = \int_{\Xi_1(\x)} \|G(\x,\xi)\|^2 \tilde{p}(\xi) d\xi, 	
\end{align} 
where the last equality follows from observing that $G(\x,\xi) = 0$ for $\xi \in \Xi_2(\x)$ and the integral in \eqref{decomp-int2} is zero because $\Xi_0(\x)$ is a measure zero set. It follows that
\begin{align} 
	& \ \mathbb{E}[\|G(\x,\xi)\|^2] \notag 
	= \int_{\Xi_1(\x)}\left(\mathcal{C}^2(2\pi\sigma^2)^{n} (4\|\xi\|_2^2 (\xi^T\x)^2) e^{-2(\xi^T\x)^2 +\us{\tfrac{\|\xi\|^2_{2}}{\sigma^2}}} \right)(2\pi\sigma^2)^{-\tfrac{n}{2}} e^{\tfrac{\us{-\|\xi\|^2_{2}}}{2\sigma^2}} d\xi  \\
	\label{ineq-lp-1}
	& \leq \int_{\Xi_1(\x)}\left(\mathcal{C}^2(2\pi\sigma^2)^{n} (4\|\xi\|_2^2 (\xi^T\x)^2) e^{-2(\xi^T\x)^2 +\us{\|\xi\|^2_{p}}} \right)(2\pi\sigma^2)^{-\tfrac{n}{2}} e^{\tfrac{\us{-\|\xi\|^2_{2}}}{2\sigma^2}} d\xi  \\
	\label{ineq-lp-2}
	& \leq \int_{\Xi_1(\x)}\left( \mathcal{C}^2(2\pi\sigma^2)^{n} (4\|\xi\|_2^2 (\xi^T\x)^2) e^{-(\xi^T\x)^2}  \right)(2\pi\sigma^2)^{-\tfrac{n}{2}} e^{\tfrac{\us{-\|\xi\|^2_{2}}}{2\sigma^2}} d\xi, 	
\end{align}
where \eqref{ineq-lp-2} follows from $\xi \in \Xi_1(\x)$ and \eqref{ineq-lp-1} follows from 
$$\frac{\|\xi\|_2^2}{\sigma^2} \leq \|\xi\|_p^2, \mbox{ where } \sigma^2 = \begin{cases}
	n^{1/2-1/p}, & p \geq 2 \\
	1. & 1 \leq p < 2
\end{cases}$$ We may then conclude that 
\begin{align}
	& \ \mathbb{E}[\|G(\x,\xi)\|^2]   \leq \int_{\Xi_1(\x)}\left(\mathcal{C}^2(2\pi\sigma^2)^{n} (\|\xi\|_2^2 (\xi^T\x)^2) e^{-(\xi^T\x)^2} \right)(2\pi\sigma^2)^{-\tfrac{n}{2}}e^{-\frac{\us{\|\xi\|^2_{2}}}{2\sigma^2}} d\xi \notag \\
	& \leq  e^{-1} \mathcal{C}^2(2\pi\sigma^2)^{n}\int_{\Xi_1(\x)}\left( \|\xi\|^2  \right)(2\pi\sigma^2)^{-\frac{n}{2}} e^{-\tfrac{\us{\|\xi\|^2_{2}}}{2\sigma^2}} d\xi \label{ineq-uni-2} \\ 
	& \leq e^{-1} \mathcal{C}^2(2\pi\sigma^2)^{n}\int_{\Real^n} \left(\|\xi\|_2^2\right)   (2\pi\sigma^2)^{-\frac{n}{2}} e^{-\tfrac{\us{\|\xi\|^2_{2}}}{2\sigma^2}} d\xi \notag
	=  e^{-1} \mathcal{C}^2(2\pi\sigma^2)^{n}\mathbb{E}[\|\xi\|^2].  
\end{align}
where~\eqref{ineq-uni-2} follows from Lemma~\ref{unimodal}.
}  \qed

\medskip

\noindent {\bf Proof of Lemma~\ref{equiv-ellipsoid}:} Suppose $(\x,\y)$ is feasible with respect to (\ref{main_prob}$_{A,{\rm ext}}^{2}$). Then $\x \in \mathcal{X}$ and is therefore feasible with (\ref{main_prob}$_{A}^{\Escr}$). In addition,
\begin{align*}
f(\x) & \triangleq \mathbb{P}\left\{ \zeta \in \Real^n \mid \zeta \in \Kscr_{\Escr}, \mid \zeta^{\intercal} \x \mid \ \le1 \right\}  =\mathbb{P} \left\{ \zeta \mid \zeta^{\intercal} U^{\intercal} \ib{\Sigma^{-1}} U \zeta \leq 1, \mid \zeta^{\intercal} \x \mid \le1 \right\}\\
& = \mathbb{P}\left\{ \zeta  \in \Real^n\mid \|\ib{\Sigma^{-1/2}} U\zeta\|^2_2 \leq 1, \mid \zeta^{\intercal} \x \mid \ \le1 \right\}\\
& = \mathbb{P} \left\{ U^{\intercal} \ib{\Sigma^{1/2}} \ib{\eta}  \in \Real^n\mid \|\ib{\eta}\|^2_2 \leq 1, \mid (U^{\intercal} \ib{\Sigma^{1/2} \eta})^{\intercal} \x \mid \ \le1 \right\}\\
& = \mathbb{P} \left\{ U^{\intercal} \ib{\Sigma^{1/2} \eta}  \in \Real^n\mid \|\ib{\eta}\|^2_2 \leq 1, \mid \ib{\eta^{\intercal} \Sigma^{1/2}} U \x \mid \ \le1 \right\}\\
& = \mathbb{P} \left\{ \ib{\eta}  \in \Real^n\mid \ib{\eta} \in \Kscr_{2}, \mid \ib{\eta}^{\intercal} \Sigma^{1/2} U \x \mid \ \le1 \right\} \triangleq \uss{g(\x)}. \end{align*}
\qed

\noindent \textbf{Proof of Proposition~\ref{bound-asym-sett_B}:}
(a) The result follows by a transformation argument.  We define a new variable ${ \tilde{\zeta} \in
\tilde{\Kscr}}$ such that ${ \tilde{\zeta}\triangleq \zeta -\mu}$ where ${
\tilde{\Kscr} \triangleq \{ \tilde{\zeta}: \| \tilde{\zeta} \|_p \leq \alpha\}}$.  The
set $\tilde{\K}(\x)$ can be defined as the following 
\ieb{
\begin{align*}
	\tilde{\K}(\x) =\left\{\tilde{\zeta}:\tilde{\zeta} \in \tilde{\Kscr} \right\}\cap \left\{\tilde{\zeta}: \tilde{\zeta} \leq T\x-\mu\right\}.
\end{align*}
}
\uss{We first show that $\zeta \in \K(\x)$ if and only if $\tilde{\zeta} \in \tilde{\K}(\x)$. Suppose $\zeta \in \K(\x)$. Then $\zeta \in \mathcal{K}$ and $c(\x,\zeta) = T\x-\zeta \geq 0$. If $\zeta \in \mathcal{K}$, then $\|\zeta -\mu\|_p \leq \alpha$ or $\|\tilde{\zeta}\|_p \leq \alpha$ where $\tilde{\zeta} = \zeta - \mu$. Furthermore, $T\x \geq \zeta$ can be rewritten as $T\x - \mu \geq \zeta - \mu$ or $T\x - \mu \geq \tilde{\zeta}$. It follows that 
$$ \tilde{\zeta} \in \tilde{\K}(\x) = \left\{ \tilde{\zeta} \mid \tilde{\zeta} \in \tilde{\Kscr} \right\} \cap \left\{ \tilde{\zeta} \mid T\x - \mu \geq \tilde{\zeta}\right\}. $$   
The reverse direction follows similarly.
}
\uss{Consequently, $\mathbb{P}\left\{ \zeta \mid \zeta \in \K(\x)\right\} = \mathbb{P}\left\{ \tilde{\zeta} \mid \tilde{\zeta} \in \tilde{\K}(\x) \right\}.$ We now proceed analyze the latter probability. }
\uss{It may be observed that}  the Minkowski functional associated with $ \tilde{\Kscr}$ is
given by $ \|\tilde{\zeta}\|_{\tilde{\Kscr}} = \tfrac{1}{\alpha}\|\tilde{\zeta}\|_p$.
\uss{Since} $ {T_{i,\bullet} \x-\mu_i \geq \delta >0 }$ for $ i=1,\hdots,d $, it follows that
\begin{align*}
\tilde{\K}(\x) &= \left\{\tilde{\zeta}: \tfrac{1}{\alpha}\|\tilde{\zeta}\|_p \leq 1 \right\} \bigcap \left\{\tilde{\zeta}: \bigcap_{i=1}^d \frac{\max\{\tilde{\zeta}_i,0\}}{T_{i,\bullet} \x-\mu_i}\leq 1  \right\}\\
&= \left\{\tilde{\zeta}: \tfrac{1}{\alpha^2}\|\tilde{\zeta}\|^{2}_p \leq 1 \right\} \bigcap \left\{\tilde{\zeta}: \bigcap_{i=1}^d \left(\frac{\max\{\tilde{\zeta}_i,0\}}{T_{i,\bullet} \x-\mu_i}\right)^2 \leq 1  \right\}\\
&= \left\{\tilde{\zeta}: \max\left\{ \frac{1}{\alpha^2}\|\tilde{\zeta}\|_{p}^2,\left(\frac{\max\{\tilde{\zeta}_1,0\}}{T_{1,\bullet}\x-\mu_1}\right)^2,\cdots, \left(\frac{\max\{\tilde{\zeta}_d,0\}}{T_{d,\bullet}\x-\mu_d}\right)^2  \right\}\leq 1 \right\}.
\end{align*}
\ieb{Since $ g_i(\x,\tilde{\zeta}) \ \triangleq  \
\left(\frac{\max\{\tilde{\zeta}_i,0\}}{T_{i,\bullet} \x-\mu_i}\right)^2 $ for
$ i=1,\ldots, d $ and $ g_{d+1}(\x, \tilde{\zeta}) \triangleq
\tfrac{1}{\alpha^2}\|\tilde{\zeta}\|^{2}_{p} $ are PHFs \uss{with} degree $ 2$, then $
g(\uss{\x}, \tilde{\zeta}) \triangleq \max \{ g_1(\x,
\tilde{\zeta}),\ldots,g_{d+1}(\x,\tilde{\zeta}) \} $ is positively homogeneous
\uss{with} degree $ 2 $. By selecting $h(\zeta) = 1$ and $\uss{\Lambda} = \tilde{\K}(\x)$, we
may invoke Lemma~\ref{las-1}, leading to the following equality. }
\begin{align} \label{eq-asym-fx-B}
f(\x) = \int_{\tilde{\K}(\x)} 1 \ d \tilde{\zeta} = \frac{1}{\mathrm{Vol}(\Kscr)} \frac{1}{\Gamma(1+d/2)} \int_{\Real^d} e^{-g(\x,\xi)} \ d\xi.
\end{align}
The equation~\eqref{eq-asym-fx-B} can be rewritten as
\begin{align*}
f(\x)&=  \int_{\mathbb{R}^d} \underbrace{\left(\mathcal{C} (2\pi\sigma^2)^{d/2} e^{-g(\x,\xi) +\frac{\|\xi\|^2_{2}}{2\sigma^2}}\right)}_{\triangleq F(\x,\xi)}\underbrace{\left( \tfrac{1}{(2\pi\sigma^2)^{d/2}} {e^{-\tfrac{\|\xi\|^2_{2}}{2\sigma^2}}}
	\right)}_{\triangleq  p(\xi)}  \ d\xi \nonumber \\ 
\notag & =   \int_{\mathbb{R}^\ieb{d}} F(\x,\xi) \  \tilde p(\xi) \ d\xi  = \mathcal{C} \ \mathbb{E}_{\tilde p(\xi)}[F(\x,\xi)], 
\mbox{ where } \mathcal{C}  \triangleq \tfrac{1}{\text{Vol}(\Kscr)} \ \ieb{\tfrac{1}{\Gamma(1+d/2)}}, 
\end{align*}
\noindent (b) Omitted (similar to proof of Lemma~\ref{bound sub settingB} (a)).

\noindent (c) When $ \Kscr$ satisfies Assumption~\ref{dist-ass-3}, the proof of Lemma~\ref{bound sub settingB}(b) requires slight modification. Suppose $ F(\x,\xi) $ and $ p(\xi) $ are defined as in (a). Then we may define $\partial F(\x,\xi)$ as 
\ieb{
\begin{align*} \scriptsize
	\partial F(\x,\xi) = \begin{cases}  \left( \ieb{\mathcal{C}}(2\pi\sigma^2)^{\ieb{d}/2} \frac{2(\max\{\xi_i,0\})^2 T_{i,\bullet}^T}{(T_{i,\bullet} \x-\mu_i)^3}  e^{-g_i(\x,\xi)+\frac{\us{\|\xi\|_{2}^2}}{2\sigma^2}}\right),  & \xi \in \Xi_i(\x), i = 1, \cdots, d \\ 
		\left( -\ieb{\mathcal{C}}(2\pi\sigma^2)^{\ieb{d}/2}  e^{-g(\x,\xi)+\frac{\us{\|\xi\|_{2}^2}}{2\sigma^2}}\right) H(\x,\xi),  & \xi \in \Xi_0(\x)\\
		{\bf 0}.  & \xi \in \Xi_{d+1}(\x), 
	\end{cases}
\end{align*}	
}
where $H(\x,\xi)$ denotes the Clarke generalized gradient of $g(\x,\xi)$, defined as in~\eqref{eq-clarke-H}.	
Consequently, it follows that $\mathbb{E} \left[\|G(\x,\xi)\|^2\right]$ is bounded as follows.
\ieb{
\begin{align}\notag
	& \quad \mathbb{E}\left[\|G(\x,\xi)\|^2\right]  = \int_{\Real^\ieb{d}} \|G(\x,\xi)\|^2 \tilde{p}(\xi) d\xi \\
	& = \sum_{i=1}^{d} \int_{\Xi_i(\x)} \|G(\x,\xi)\|^2 \tilde{p}(\xi) d\xi \notag 	
	+ \int_{\Xi_{d+1}(\x)} \|\underbrace{G(\x,\xi)}_{ \ = \ {\bf 0}}\|^2 \tilde{p}(\xi) d\xi 	
	\\ & + \int_{\Xi_0(\x)} \|G(\x,\xi)\|^2 \tilde{p}(\xi) d\xi \label{decomp-int4}\\
	\notag			&  = \sum_{i=1}^{d}\int_{\Xi_i(\x)} \|G(\x,\xi)\|^2 \tilde{p}(\xi) d\xi, 	
\end{align} 
where the last equality follows from observing that $G(\x,\xi) = 0$ for $\xi \in \Xi_{d+1}(\x)$ and the integral in \eqref{decomp-int4} is zero because $\Xi_0(\x)$ is a measure zero set. }
It follows that
\begin{align*}
\mathbb{E}&[\|G(x,\xi)\|^2] \\
&=\sum_{i=1}^{d}\int_{\Xi_i(x)} 4\mathcal{C}^2(2\pi\sigma^2)^d \frac{\|T_{i,\bullet}\|^2}{(T_{i,\bullet} \x-\mu_i)^2}\left(\frac{\xi_k}{{T_{i,\bullet} \x-\mu_i}}\right)^4 e^{-\frac{2(\xi_i)^2}{(T_{i,\bullet} \x-\mu_i)^2}+\frac{\|\xi\|_{2}^2}{\sigma^2}}\left( \tfrac{1}{(2\pi\sigma^2)^{d/2}} e^{-\tfrac{\|\xi\|_{2}^2}{2\sigma^2}}\right)d\xi,\\
&\leq \sum_{i=1}^{d}\int_{\Xi_i(x)} 4\mathcal{C}^2(2\pi\sigma^2)^d \frac{\|T_{i,\bullet}\|^2}{\delta^2}\left(\frac{\xi_k}{{T_{i,\bullet} \x-\mu_i}}\right)^4 e^{-\frac{2(\xi_i)^2}{(T_{i,\bullet} \x-\mu_i)^2}+\frac{\|\xi\|_{2}^2}{\sigma^2}}\left( \tfrac{1}{(2\pi\sigma^2)^{d/2}} e^{-\tfrac{\|\xi\|_{2}^2}{2\sigma^2}}\right)d\xi\\
&\leq \sum_{i=1}^{d}\int_{\Xi_i(x)} 4\mathcal{C}^2(2\pi\sigma^2)^d \frac{\|T_{i,\bullet}\|^2}{\delta^2}\left(\frac{\xi_i}{{T_{i,\bullet} \x-\mu_i}}\right)^4 e^{-\left(2-\tfrac{\alpha^2}{\sigma^2}\right)\frac{(\xi_i)^2}{(T_{i,\bullet} \x-\mu_i)^2}}\left( \tfrac{1}{(2\pi\sigma^2)^{d/2}} e^{-\tfrac{\|\xi\|_{2}^2}{2\sigma^2}}\right)d\xi
\end{align*}
where the first inequality follows from $ {T_{i,\bullet} \x-\mu_i \geq \delta >0 }$ for all $ i$, and the second inequality follows from $ \xi \in \Xi_i(\x)$. 
It follows from Lemma~\ref{unimodal} that given any $ \alpha$, choosing the variance $ \sigma^2$ of $ \xi $ such that $\sigma^2 = \alpha^2 $ leads to the bound
$
{\mathbb{E} [\|G(\x,\xi)\|^2] \leq 16\mathcal{C}^2(2\pi\sigma^2)^d \uss{\sum_{i=1}^d} \frac{\|T_{i,\bullet}\|^2}{\delta^2 e^2}}.
$ \qed

\medskip

\noindent \textbf{Proof of Lemma~\ref{cond-sec-mom}:}
\noindent	{If $\tilde G(\x_k,\xi) \triangleq G(\x_k,\xi) - \mathbb{E}[G(\x_k,\xi)]$, by the conditional independence of $\tilde{G}(\x_k,\xi_j)$ and $\tilde{G}(\x_k,\xi_{\ell})$ for $j \neq \ell$, we have 
\begin{align} \notag
	& \quad \mathbb{E}[\|\bar{w}_{G,k}\|^2 \mid \mathcal{F}_k ] = \frac{1}{N^2_k}\mathbb{E}\left[\left\|\sum_{j=1}^{N_k} \tilde{G}(\x_k,\xi_j)\right\|^2 \mid \mathcal{F}_k \right] \\
	&  \notag= \frac{1}{N_k^2} \mathbb{E}\left[\left[\sum_{j=1}^{N_k} \|\tilde{G}(\x_k,\xi_j)\|^2 + \sum_{\ell \neq j} 2\tilde{G}(\x_k,\xi_{\ell})^T \tilde{G}(\x_k,\xi_j) \right]\mid \mathcal{F}_k \right] \\ 
	& \notag = \frac{1}{N_k}\left( \mathbb{E}\left[ \|{G}(\x_k,\xi)\|^2\mid \mathcal{F}_k \right] + \| \mathbb{E}[G(\x_k,\xi) \mid \mathcal{F}_k]\|^2 - 2\mathbb{E}\left[G(\x_k,\xi)\mid \mathcal{F}_k\right]^T\mathbb{E}[G(\x_k,\xi) \mid \mathcal{F}_k]\right)  \\
	& \label{def-wg} \leq \frac{1}{N_k} \mathbb{E} \left[\|{G}(\x_k,\xi)\|^2 \mid \mathcal{F}_k\right].
\end{align}
{By \eqref{def-wg} \uss{and Prop.~\ref{bound sub}}, 
	$\mathbb{E}[\|\bar{w}_{G,k}\|^2 \mid \mathcal{F}_k ]  \uss{ \ \leq \ }	 \frac{\mathcal{C}^2_{\Kscr}(2\pi)^n}{eN_k} \mathbb{E}_{\tilde p}[\|\xi\|^2]$ \uss{for Setting A.}
	Similarly, \uss{for Setting B}, \uss{by Lemma.~\ref{bound sub settingB}},  
	$$\uss{\mathbb{E}[\|\bar{w}_{G,k}\|^2 \mid \mathcal{F}_k ] 
		\leq	16\mathcal{C}^2(2\pi\sigma^2)^d \uss{\sum_{i=1}^d} \frac{\|T_{i,\bullet}\|^2}{\delta^2 e^2 N_k} }.$$ In addition, for Setting A, 
	$ \mathbb{E}[\|\bar{w}_{f,k}\|^2 \mid \mathcal{F}_k ] 
	\leq \frac{2(\mathcal{C}_{\Kscr}^2(2\pi)^n+1)}{N_k} \mbox{ and } 
	\mathbb{E}[\|\bar{w}_{f,k}\|^2 \mid \mathcal{F}_k ] 
	\leq \frac{\uss{\mathcal{C}^2(2\pi\sigma^2)^d}}{N_k}.$ \qed 
}

\medskip

\noindent \textbf{Proof of Lemma~\ref{bd-sec-moment}:}
(Setting A) Consider $\bar{w}_k$, defined as
$ \bar{w}_k \triangleq \frac{\ib{-}(G_k+\bar{w}_{G,k})}{(f(\x_k) + \bar{w}_{f,k})^{2}+\epsilon_k}- \frac{\ib{-}G_k}{(f(\x_k))^{2}}.$
We have that 
\begin{align*}
\|\bar{w}_k\|^2 & = \left\|\frac{\ib{-}(G_k+\bar{w}_{G,k})}{(f(\x_k) + \bar{w}_{f,k})^{2}+\epsilon_k}- \frac{\ib{-}G_k}{(f(\x_k))^{2}} \right\|^2 \\
& =  \left\|\frac{\ib{-}(G_k+\bar{w}_{G,k})}{(f(\x_k) + \bar{w}_{f,k})^{2}+\epsilon_k}-\frac{\ib{-}(G_k+\bar{w}_{G,k})}{(f(\x_k))^{2}+\epsilon_k}+\frac{\ib{-}(G_k+\bar{w}_{G,k})}{(f(\x_k))^{2}+\epsilon_k}\right. \\ 
& \left.- \frac{\ib{-}G_k}{(f(\x_k))^{2}+\epsilon_k}
+ \frac{\ib{-}G_k}{(f(\x_k))^{2}+\epsilon_k}- \frac{\ib{-}G_k}{(f(\x_k))^{2}} \right\|^2 \\
& \leq  3\left\|G_k - G_k+\bar{w}_{G,k}\right\|^2 \frac{1}{((f(\x_k))^{2}+\epsilon_k)^2}\\
&+3\left\|G_k+\bar{w}_{G,k}\right\|^2 \left\|\frac{1}{(f(\x_k))^{2}+\epsilon_k} -\frac{1}{(f(\x_k) + \bar{w}_{f,k})^{2}+\epsilon_k}\right\|^2 \\
&+ 3\left\|G_k\right\|^2 \left\|\frac{1}{(f(\x_k))^{2}} -\frac{1}{(f(\x_k))^{2}+\epsilon_k}\right\|^2 \\
& \leq  3\left\|\bar{w}_{G,k}\right\|^2 \frac{1}{((f(\x_k))^{2}+\epsilon_k)^2}\\
&+3\left\|G_k+\bar{w}_{G,k}\right\|^2 \left\|\frac{(2f(\x_k)+\bar{w}_{f,k}) \bar{w}_{f,k}}{((f(\x_k))^{2}+\epsilon_k)((f(\x_k) + \bar{w}_{f,k})^{2}+\epsilon_k)}\right\|^2 \\
& + 3\left\|G_k\right\|^2 \left\|\frac{\epsilon_k}{(f(\x_k))^{2}((f(\x_k))^{2}+\epsilon_k)}\right\|^2 \\
& \leq  3\left\|\bar{w}_{G,k}\right\|^2 \frac{1}{\epsilon_f^{{4}}}+{3}\left\|G_k+\bar{w}_{G,k}\right\|^2 \left\|\frac{(2f(\x_k)+\bar{w}_{f,k})}{\epsilon_f^{2} \epsilon_k}\right\|^2\|\bar{w}_{f,k}\|^2 
+ 3\left\|G_k\right\|^2 \left(\frac{\epsilon_k^2}{\epsilon_f^{{8}}}\right) \\
& \leq  3\left\|\bar{w}_{G,k}\right\|^2 \frac{1}{\epsilon_f^4}+{3} \left\|G_k+\bar{w}_{G,k}\right\|^2 \frac{(8f^2(\x_k)\|\bar{w}_{f,k}\|^2+2\|\bar{w}_{f,k}\|^4)}{\epsilon_f^{4} \epsilon^2_k}   + \left\|G_k\right\|^2 \left(\frac{\epsilon_k^2}{\epsilon_f^{8}}\right), 
\end{align*}
where $f(\x_k) \geq \epsilon_f$ for every $\x_k \in \mathcal{X}.$ Taking conditional expectations  and recalling the independence of $\bar{w}_{f,k}$ and $\bar{w}_{G,k}$ conditional on $\mathcal{F}_k$, the following bound emerges. 
\begin{align*}
\mathbb{E}&[\|\bar{w}_k\|^2 \mid \mathcal{F}_k] 
\leq 3\mathbb{E}\left[\left\|\bar{w}_{G,k}\right\|^2 \mid \mathcal{F}_k\right] \frac{1}{\epsilon_f^2}  \\
&+{3}\mathbb{E}\left[\left\|G_k+\bar{w}_{G,k}\right\|^2 \frac{(8f^2(\x_k)\|\bar{w}_{f,k}\|^2+2\|\bar{w}_{f,k}\|^4)}{\epsilon_f^{{4}} \epsilon_k^2}\mid \mathcal{F}_k\right] 
+ 3\mathbb{E}\left[\left\|G_k\right\|^2 \left(\frac{\epsilon_k^2}{\epsilon_f^{{8}}}\right)\mid \mathcal{F}_k\right] \\
& \leq 	 3\frac{\nu_G^2}{\epsilon_f^2N_k} +{3}\mathbb{E}\left[\left\|G_k+\bar{w}_{G,k}\right\|^2 \mid \mathcal{F}_k\right]  \mathbb{E}\left[ \frac{(8f^2(\x_k)\|\bar{w}_{f,k}\|^2+2\|\bar{w}_{f,k}\|^4)}{\epsilon_f^{{4}} \epsilon_k^2}\mid \mathcal{F}_k\right]
+ \left(\frac{3\epsilon_k^2M_G^2}{\epsilon_f^{{8}}}\right) \\
& \leq 	 3\frac{\nu_G^2}{\epsilon_f^2N_k}+{3} M_G^2 \frac{8f^2(\x_k)\nu_f^2}{\epsilon_f^{4} \epsilon_k^2N_k}+{3}M_G^2\mathbb{E}\left[\frac{\|\bar{w}_{f,k}\|^4}{\epsilon_f^{4} \epsilon_k^2}\mid \mathcal{F}_k\right] +
3\left(\frac{\epsilon_k^2M_G^2}{\epsilon_f^{{8}}}\right), 
\end{align*} 
where {$\|G_k\|^2 = \| \mathbb{E}[G(\x_k,\xi) \mid \mathcal{F}_k]\|^2 \leq \mathbb{E}[\| G(\x_k,\xi)\|^2 \mid \mathcal{F}_k] \leq M_G^2$ by Jensen's inequality.}
From~Prop.~\ref{bound sub}(b,c), $| F(\x,\xi)| \leq \ib{M_F}$ for any $\x, \xi$, {implying that}
\begin{align*}
\|\bar{w}_{f,k}\|^2 &= \left\| \frac{\sum_{j=1}^{N_k} F(\x_k,\xi_j)}{N_k} - f(\x_k)\right\|^2 
\leq 2 \left\|\frac{\sum_{j=1}^{N_k} F(\x_k,\xi_j)}{N_k}\right\|^2 + 2f^2(\x_k) \leq 2(\ib{M_F}^2+1).
\end{align*}
Consequently,  by recalling that $\epsilon_k = 1/N_k^{1/4}$, the following holds a.s. 
\begin{align*}
\mathbb{E}[\|\bar{w}_k\|^2 \mid \mathcal{F}_k] & \leq \frac{3\nu_G^2}{\epsilon_f^2N_k}+{24}M_G^2   \frac{f^2(\x_k)\nu_f^2}{\epsilon_f^{{4}} \epsilon_k^2N_k}
+{3}\mathbb{E}\left[\frac{\|\bar{w}_{f,k}\|^4}{\epsilon_f^4 \epsilon_k^2}\mid \x_k\right] 
+ \left(\frac{\epsilon_k^2M_G^2}{\epsilon_f^8}\right) \\ 
& \leq \frac{\nu_G^2}{\epsilon_f^2N_k}+{24}M_G^2  \frac{f^2(\x_k)\nu_f^2}{\epsilon_f^{{4}} \epsilon_k^2N_k}+\frac{{6}({M_F}^2+1)M_G^2\nu_f^2}{\epsilon_f^4 \epsilon_k^2 N_k}\\
& \leq \frac{3\nu_G^2}{\epsilon_f^2\sqrt{N_k}}+M_G^2  \frac{{24}f^2(\x_k)\nu_f^2}{\epsilon_f^{{4}} \sqrt{N_k}}+\frac{{6}(M_F^2+1)\nu_f^2}{\epsilon_f^{{4}} \sqrt{N_k}}
+ \left(\frac{3M_G^2}{\epsilon_f^{{8}}\sqrt{N_k}}\right) \\
& \triangleq \frac{\nu^2}{\sqrt{N_k}}, \mbox{ where } 
\nu^2 \triangleq \frac{3\nu_G^2}{\epsilon_f^2}+M_G^2  \frac{{24}\nu_f^2}{\epsilon_f^{{4}} }+\frac{{6}({M_F}^2+1)\nu_f^2}{\epsilon_f^{{4}}}
+ \left(\frac{3M_G^2}{\epsilon_f^{{8}}}\right).
\end{align*}
\noindent (Setting B)  \uss{Since  
$ \bar{w}_k  \triangleq {\frac{\ib{-}(G_{k}+\bar{w}_{G,k})}{(f(\x_k) + \bar{w}_{f,k})+\epsilon_k}+ \frac{G_{k}}{f(\x_k)}}$ and 
\begin{align*}
	\|\bar{w}_{k}\|^2 & = \left\|\frac{\ib{-}(G_{k}+\bar{w}_{G,{k}})}{(f(\x_k) + \bar{w}_{f,{k}})+\epsilon_k}- \frac{\ib{-}G_k}{(f(\x_k))} \right\|^2 \\
	& =  \left\|\frac{\ib{-}(G_{k}+\bar{w}_{G,{k}})}{(f(\x_k) + \bar{w}_{f,k})+\epsilon_k}-\frac{\ib{-}(G_{k}+\bar{w}_{G,{k}})}{f(\x_k)+\epsilon_k}+\frac{\ib{-}(G_{k}+\bar{w}_{G,k})}{f(\x_k)+\epsilon_k}\right. \\ 
	& \left.- \frac{\ib{-}G_{k}}{f(\x_k)+\epsilon_k}
	+ \frac{\ib{-}G_{k}}{f(\x_k)+\epsilon_k}- \frac{\ib{-}G_{k}}{f(\x_k)} \right\|^2 \\
	& \leq  3\left\|G_{k} - G_{k}+\bar{w}_{G,{k}}\right\|^2 \frac{1}{{(f(\x_k)+\epsilon_k)^2}}
	\\
	& +3\left\|G_{k}+\bar{w}_{G,{k}}\right\|^2 \left\|\frac{1}{f(\x_k)+\epsilon_k} -\frac{1}{(f(\x_k) + \bar{w}_{f,{k}})+\epsilon_k}\right\|^2 \\
	&+ 3\left\|G_{k}\right\|^2 \left\|\frac{1}{f(\x_k)} -\frac{1}{f(\x_k)+\epsilon_k}\right\|^2 \\
	& \leq  3\left\|\bar{w}_{G,{k}}\right\|^2 \frac{1}{{(f(\x_k)+\epsilon_k)^2}}
	+3\left\|G_{k}+ \bar{w}_{f,{k}}\right\|^2 \left\| \frac{\bar{w}_{f,k}}{(f(\x_k)+\epsilon_k)(\underbrace{(f(\x_k) + \bar{w}_{f,{k}})}_{{\geq 0, F(\x_k,\xi) \geq 0}}+\epsilon_k)}\right\|^2 \\
	& + 3\left\|G_{k}\right\|^2 \left\|\frac{\epsilon_k}{f(\x_k)(f(\x_k)+\epsilon_k)}\right\|^2 \\
	& \leq  3\left\|\bar{w}_{G,{k}}\right\|^2 \frac{1}{{\epsilon_f^2}}
	+3\left\|G_k+\bar{w}_{G,{k}}\right\|^2 \left\|\frac{1}{\epsilon_f \epsilon_k}\right\|^2\|\bar{w}_{f,{k}}\|^2 
	+ 3\left\|G_{k}\right\|^2 \left(\frac{\epsilon_k^2}{\epsilon_f^{4}}\right) \\
	& \leq  3\left\|\bar{w}_{G,{k}}\right\|^2 \frac{1}{\epsilon_f^2}+3\left\|G_{k}+\bar{w}_{G,{k}}\right\|^2 \frac{\|\bar{w}_{f,{k}}\|^2}{\epsilon_f^{2} \epsilon^2_k} + \left\|G_{k}\right\|^2 \left(\frac{\epsilon_k^2}{\epsilon_f^{4}}\right), 
\end{align*}
where $f(\x_k) \geq \epsilon_f$ and for every $\x_k \in \mathcal{X}.$ Taking expectations conditioned on $\mathcal{F}_k$ and recalling the independence of $\bar{w}_{f,k}$ and $\bar{w}_{G,k}$ conditional on $\mathcal{F}_k$, we have the following bound. 
\begin{align*}
	\mathbb{E}&[\|\bar{w}_{k}\|^2 \mid \mathcal{F}_k] \\
	& \leq \left(3\mathbb{E}\left[\left\|\bar{w}_{G,{k}}\right\|^2 \mid \mathcal{F}_k\right] \frac{1}{\epsilon_f^2}  +3\mathbb{E}\left[\left\|G_{k}+\bar{w}_{G,{k}}\right\|^2  \frac{\|\bar{w}_{f,{k}}\|^2}{\epsilon_f^{2} \epsilon_k^2}\mid \mathcal{F}_k\right] 
	+ 3\mathbb{E}\left[\left\|G_{k}\right\|^2 \left(\frac{\epsilon_k^2}{\epsilon_f^{4}}\right)\mid \mathcal{F}_k\right]\right) \\
	& \leq 	\left( 3\frac{\nu_{G}^2}{\epsilon_f^2N_k}  
	+ 3\mathbb{E}\left[\left\|G_{k}+\bar{w}_{G,{k}}\right\|^2 \mid \mathcal{F}_k\right] \mathbb{E}\left[ \frac{\|\bar{w}_{f,{k}}\|^2}{\epsilon_f^{2} \epsilon_k^2}\mid \mathcal{F}_k\right] 
	+ 3\left(\frac{\epsilon_k^2M_{G}^2}{\epsilon_f^{4}}\right) \right)\\
	& \leq 	  \left(3\frac{\nu_{G}^2}{\epsilon_f^2N_k}+3M_G^2  \frac{\nu_{f}^2}{\epsilon_f^{2} \epsilon_k^2N_k}+ 
	3\left(\frac{\epsilon_k^2M_{G}^2}{\epsilon_f^{4}}\right)\right). 
\end{align*} 
By selecting  $\epsilon_k = 1/N_k^{1/4}$, we have that
\begin{align*}
	\mathbb{E}[\|\bar{w}_{k}\|^2 \mid \mathcal{F}_k] & \leq  \frac{\nu^2}{\sqrt{N_k}}, \mbox{ where } 
	\nu^2 \triangleq \left(3\frac{\nu_{G}^2}{\epsilon_f^2}+3M_G^2  \frac{\nu_{f}^2}{\epsilon_f^{2} }+ 
	3\left(\frac{M_{G}^2}{\epsilon_f^{4}}\right)\right).
\end{align*}}
\noindent \textbf{Proof of Proposition~\ref{bd-itercomp-sc-sgd}:}
(i) Using the update rule of $ \x_{k+1}$ and the fact that $\x^*=\Pi_{\mathcal X} [\x^*]$, for any $d_k + \bar{w}_{k}$ where $d_k \in \partial h(\x_k)$ and $k \geq 1$,
\begin{align*}
{1\over 2}\|\x_{k+1}-\x^*\|^2&\us{ = } {1\over 2}\|\Pi_{\mathcal X}(\x_k-\gamma_k \us{(d_k + \bar{w}_k)})-\Pi_{\mathcal{X}}(\x^*))\|^2\\
&\leq {1\over 2}\|\x_k-\gamma_k \us{(d_k+\bar{w}_k)}-\x^*\|^2\\
&={1\over 2}\|\x_k-\x^*\|^2+{1\over 2}\gamma_k^2\|\us{d_k + \bar{w}_k}\|^2-\gamma_k(\x_k-\x^*)^T(\us{d_k}+\bar w_{k}),
\end{align*}
where in the second inequality, we {employ the non-expansivity of projection operator}. Now by using the convexity of $h$, we obtain:
\begin{align*}
2\gamma_k (h(\x_k)-h(\x^*))&\leq \left(\|\x_k-\x^*\|^2-\|\x_{k+1}-\x^*\|^2\right)+{\|\us{d_k+\bar{w}_k}\|^2\gamma_k^2}\\
&-{2}\gamma_k\bar w_{k}^T(\x_k-\x^*)\\
& \leq \left(\|\x_k-\x^*\|^2-\|\x_{k+1}-\x^*\|^2\right)+{\|\us{d_k+\bar{w}_k}\|^2\gamma_k^2}\\
&+ \gamma_k^2 \|\x_k-\x^*\|^2 + \|\bar{w}_k\|^2,
\end{align*}
where we use $a^Tb\leq{1\over 2}\|a\|^2+{1\over 2}\|b\|^2$. Now by summing from $k = \widehat{K}$ to $K-1$, where $\widehat{K}$ is an integer satisfying $0 \leq \widehat{K} < K-1$,  we obtain the next inequality.
\begin{align*}
\sum_{k=\widehat{K}}^{K-1} 2\gamma_k (h(\x_k)-h(\x^*))
& \leq {\|\x_{\hat K}-\x^*\|^2}+\sum_{k=\widehat{K}}^{K-1} \gamma_k^2({\|\us{d_k+\bar{w}_k}\|^2}+ \|\x_k-\x^*\|^2)  +\|\bar{w}_k\|^2.
\end{align*}
Dividing both sides by $2\sum_{k=\widehat{K}}^{K-1} \gamma_k$, taking expectations on both sides,  
{and invoking Lemma~\ref{bd-sec-moment} which leads to}  $\mathbb E[\|\bar w_{k}\mid \mathcal F_k\|]^2\leq {\nu^2\over \sqrt{N_k}}$ and the bound of the subgradient, i.e., \us{$\mathbb E[\|d_k + \bar{w}_k\|^2]\leq M_G^2$}, we obtain the following bound.
\begin{align}\notag
\mathbb{E}&\left[\frac{\sum_{k=\widehat{K}}^{K-1} 2\gamma_k (h(\x_k)-h(\x^*))}{\sum_{k=\widehat{K}}^{K-1} 2\gamma_k}\right]\\
& \leq \mathbb{E}\left[\frac{{\|\x_{\hat K}-\x^*\|^2}+\sum_{k=\widehat{K}}^{K-1} \gamma_k^2{\|\us{d_k+\bar{w}_k}\|^2}+ \sum_{k=\widehat{K}}^{K-1} \gamma_k^2 \|\x_k-\x^*\|^2 + \sum_{k=\widehat{K}}^{K-1} \|\bar{w}_k\|^2}{\sum_{k=\widehat{K}}^{K-1}2\gamma_k}\right] \\
& \leq \frac{\mathbb{E}[\|x_{\widehat{K}}-x^*\|^2]}{\sum_{k=\widehat{K}}^{K-1}2\gamma_k } + \frac{\sum_{k=\widehat{K}}^{K-1}\gamma_k^2 (M_G^2 + B^2)}{\sum_{k=\widehat{K}}^{K-1}2\gamma_k }+ \frac{\sum_{k=\widehat{K}}^{K-1} \tfrac{\nu^2}{\sqrt{N_k}}}{\sum_{k=\widehat{K}}^{K-1}2\gamma_k}. \label{bd-prop} 
\end{align}
By utilizing Jensen's inequality, we obtain that 
$$  \mathbb{E}\left[(h(\bar{x}_{\widehat{K},K}-h(\x^*))\right] \leq \mathbb{E}\left[\frac{\sum_{k=\widehat{K}}^{K-1} 2\gamma_k (h(\x_k)-h(\x^*))}{\sum_{k=\widehat{K}}^{K-1} 2\gamma_k}\right],$$ 
where $\bar{x}_{\widehat{K},K} \triangleq \tfrac{\sum_{k=\widehat{K}}^{K-1} \gamma_k x_k}{\sum_{k=\widehat{K}}^{K-1} \gamma_k},$ which when combined with \eqref{bd-prop} leads to  \eqref{bd-prop-h}. \qed

\bibliographystyle{spmpsci}      
\bibliography{mybib}
\end{document}